\documentclass{amsart}
\usepackage{graphicx}
\usepackage[all]{xy}
\usepackage{latexsym}
\usepackage{amssymb}
\usepackage{amsmath}
\usepackage{color}

\usepackage{fullpage}
\xyoption{arc}
 \newfam\cyrfam

  \font\tencyr=wncyr10

  \font\sevencyr=wncyr7

  \font\fivecyr=wncyr5

  \textfont\cyrfam=\tencyr \scriptfont\cyrfam=\sevencyr

    \scriptscriptfont\cyrfam=\fivecyr

  \newfam\cyifam

  \font\tencyi=wncyi10

  \font\sevencyi=wncyi7

  \font\fivecyi=wncyi5

  \textfont\cyifam=\tencyi \scriptfont\cyifam=\sevencyi

    \scriptscriptfont\cyifam=\fivecyi

\def\id{{\mbox{1 \hskip -7pt 1}}}
\newcommand{\sgn}{{\mathit s  \mathit g\mathit  n}}
 \newcommand{\lon}{\longrightarrow}
 \newcommand{\bu}{\bullet}
 
 \newcommand{\rar}{\rightarrow}

 \newcommand{\Z}{{\mathbb Z}}
 \newcommand{\bS}{{\mathbb S}}

 \newcommand{\R}{{\mathbb R}}
 \newcommand{\N}{{\mathbb N}}
 \newcommand{\K}{{\mathbb K}}

 \newcommand{\ot}{\otimes}
 
\newcommand{\GC}{\mathsf{GC}}
\newcommand{\OGC}{\mathsf{OGC}}
\newcommand{\dGC}{\mathsf{dGC}}

 \newcommand{\Beq}{\begin{equation}}
 \newcommand{\Eeq}{\end{equation}}
 \newcommand{\Beqr}{\begin{eqnarray}}
 \newcommand{\Eeqr}{\end{eqnarray}}
 \newcommand{\Beqrn}{\begin{eqnarray*}}
 \newcommand{\Eeqrn}{\end{eqnarray*}}
 \newcommand{\Ba}{\begin{array}}
 \newcommand{\Ea}{\end{array}}
 \newcommand{\Bi}{\begin{itemize}}
 \newcommand{\Ei}{\end{itemize}}
 \newcommand{\Bc}{\begin{center}}
 \newcommand{\Ec}{\end{center}}
 \newcommand{\cA}{{\mathcal A}}
 \newcommand{\cB}{{\mathcal B}}
 
 \newcommand{\caD}{{\mathcal D}}
 \newcommand{\cE}{{\mathcal E}}

 \newcommand{\cI}{{\mathcal I}}

 \newcommand{\caL}{{\mathcal L}}
 \newcommand{\cM}{{\mathcal M}}

 \newcommand{\ga}{\gamma}
 
 \newcommand{\Ga}{\Gamma}

 \newcommand{\sip}{\smallskip}
 \newcommand{\bip}{\bigskip}
 \newcommand{\mip}{\vspace{2.5mm}}
\newcommand{\RGC}{\mathsf{RGC}}

\newcommand{\Lie}{\mathcal{L} \mathit{ie}}

 \newcommand{\bbu}{\mbox{\resizebox{2.8mm}{!}{$\bullet$}}}

 \newcommand{\Ass}{\mathcal{A}ss}

\newcommand{\ORGC}{\mathsf{ORGC}}
\newcommand{\dRGC}{\mathsf{dRGC}}
\newcommand{\sRGC}{\Delta\mathsf{RGC}}

\theoremstyle{plain}
\swapnumbers

\newtheorem{prop-def}[theorem]{Proposition-definition}

\newtheorem{f-theorem}{Formality Theorem}[section]
\newtheorem{main-theorem}{Main~Theorem}[section]
\newtheorem{section-theorem}{Theorem}[section]

\theoremstyle{definition}

\usepackage{tikz}
\usetikzlibrary{matrix}
\usetikzlibrary{shapes}
\usetikzlibrary{arrows}
\usetikzlibrary{calc,3d}
\usetikzlibrary{decorations,decorations.pathmorphing,decorations.pathreplacing,decorations.markings}
\usetikzlibrary{through}

\tikzset{ext/.style={circle, draw,inner sep=1pt},int/.style={circle,draw,fill,inner sep=1.4pt},nil/.style={inner sep=1pt}}
\tikzset{cy/.style={circle,draw,fill,inner sep=2pt},scy/.style={circle,draw,inner sep=2pt},scyx/.style={draw,cross out,inner sep=2pt},scyt/.style={draw,regular polygon,regular polygon sides=3,inner sep=0.95pt}}
\tikzset{exte/.style={circle, draw,inner sep=3pt},inte/.style={circle,draw,fill,inner sep=3pt}}
\tikzset{diagram/.style={matrix of math nodes, row sep=3em, column sep=2.5em, text height=1.5ex, text depth=0.25ex}}
\tikzset{diagram2/.style={matrix of math nodes, row sep=0.5em, column sep=0.5em, text height=1.5ex, text depth=0.25ex}}

 \tikzset{
  rightblack/.style={
    decoration={markings,mark=at position .8 with {\arrow[scale=1.2,black]{latex}}},
    postaction={decorate},
    shorten >=0.4pt}}
\tikzset{
  leftblack/.style={
    decoration={markings,mark=at position .55 with {\arrowreversed[scale=1.2,black]{latex}}},
    postaction={decorate},
    shorten >=0.4pt}}

\begin{document}

 \sloppy

 \newenvironment{proo}{\begin{trivlist} \item{\sc {Proof.}}}
  {\hfill $\square$ \end{trivlist}}

\long\def\symbolfootnote[#1]#2{\begingroup%
\def\thefootnote{\fnsymbol{footnote}}\footnote[#1]{#2}\endgroup}

\title{A low-valence ribbon graph complex\\ computing the cohomology of $\cM_{g,m}$}

\author{Sergei A.\ Merkulov}
\address{Sergei~A.\ Merkulov,  Faculty of Mathematics, National Research University Higher School of Economics, Moscow }
\email{smerkulov25@gmail.com}

\begin{abstract} It is proven that every cohomology class of the moduli space $\cM_{g,m}$ for any $2g+m\geq 3$, $m\geq 1$ can be represented combinatorially by a ribbon quiver with at most four-valent vertices. The {\it at most four}-valency condition is sharp.

%\bip

%\noindent {\sc Mathematics Subject Classifications} (2000). 14H10, 18G85, %18M70

%\noindent {\sc Key words}. Moduli spaces of algebraic curves, graph complexes.
\end{abstract}

 \maketitle
\markboth{}{}

{\Large
\section{\bf Introduction}
}

\subsection{Graph complexes} The famous graph complex $\GC_d$, $d\in \Z$, was introduced in \cite{Ko1} in the context 
of deformation quantizations of Poisson structures and found applications in many branches of mathematics (see \cite{Ko2,Wi3} for overviews). It is generated by graphs with at least trivalent vertices (which are assigned the cohomological degree $d$) and with edges  (which are assigned the cohomological degree $1-d$) whose directions can be flipped with the sign factor $(-1)^d$ for each flip. There are essentially two complexes in the family $\{\GC_d\}_{d\in \Z}$ with non-isomorphic cohomologies, one for $d$ even and one for $d$ odd. 

\sip

Fixing a direction on each edge and assuming that valencies of vertices are at least 2 (and at least one vertex has valency $\geq 3$) gives us a directed version $\dGC_d$ of $\GC_d$.  There is a quasi-isomorphism of complexes \cite{Wi1}
\Beq\label{1: GC to dGC}
\GC_d\lon \dGC_d
\Eeq
which sends a graph $\Ga\in \GC_d$ with no fixed directions on edges into a sum of graphs obtained from $\Ga$ by choosing directions on edges in all possible ways. Thus the complex $\dGC_d$ gives us essentially nothing new, but it contains an interesting subcomplex $\OGC_d\subset \dGC_d$ spanned by graphs with no closed paths of directed edges (such graphs are called often {\it quivers} or {\it oriented graphs}). It was proven in \cite{Wi2,Z,MWW} that
there is a $\Lie_\infty$ quasi-isomorphism of the two dg Lie algebras (note the shift in $d$),
$$
\GC_d \simeq \OGC_{d+1}.
$$
The complexes $\OGC_{d}$ controls the homotopy theory of (degree shifted)
Lie bialgebras \cite{MW1} and  their deformation quantizations \cite{MW2}; they are also used in the study of the real locus of the moduli space of genus zero curves [KW]. The complex $\OGC_{d+1}$ (and hence $\GC_d$) admits a surprisingly  ``small" model $\GC_d^{\mathsf T}$ which is spanned by graphs whose all vertices are precisely   {\it trivalent} (see \S 3.3 in \cite{Me1}).

\subsection{Ribbon graph complexes}
There is a cochain complex $\RGC_0$  spanned by  ribbon graphs with at least trivalent vertices and marked boundaries which computes the total compactly supported cohomology  \cite{Pe, Ko0}
\Beq\label{1: RGC_0 versus H(M)}
H^\bu(\RGC_0)=\prod_{g\geq 0, m\geq 1\atop 2g+m\geq 3} H^{\bu}_c(\cM_{g,m} \times (\R^{>0})^m)\simeq \prod_{g\geq 0, m\geq 1\atop 2g+m\geq 3} H^{\bu -m}_c(\cM_{g,m})
\Eeq
of the moduli spaces $\cM_{g,m}$ of genus $g$ algebraic curves with $m$ marked points. By analogy to $\GC_d$, it is useful for applications (see e.g. \S 4 in \cite{MW0}) to consider Penner's original complex $\RGC_0$ as a member of a family of ribbon graph complexes $\RGC_d$ parameterized by integers $d\in \Z$ such that vertices
of the generating  ribbon graphs have degree $d$,  the edges are directed (up to a flip and multiplication by $(-1)^d$) and have  degree $1-d$. Contrary to the case of the Kontsevich graph complexes $\GC_d$,  {\it all}\, the complexes  $\RGC_d$ have isomorphic cohomology \cite{MW0},
$$
H^\bu(\RGC_d)=\prod_{g\geq 0, m\geq 1\atop 2g+m\geq 3} H^{\bu +(d-1)m+d(2g-1)}_c(\cM_{g,m}).
$$
Again, as in the case of $\GC_d$ one can consider a version $\mathsf{dRGC}_d$ of Penner's complex $\RGC_d$ spanned by ribbon graphs whose edges have a fixed direction and valencies of vertices are at least 2 (and at least one vertex has valency $\geq 3$). This directed complex gives us nothing new as there is a quasi-isomorphism
\Beq\label{1: RGC to dRGC}
\RGC_d \stackrel{\simeq}\lon \mathsf{dRGC}_d,
\Eeq
which is very similar to the one in the  ``commutative" case (\ref{1: GC to dGC}). The complex
$\mathsf{dRGC}_d$ contains a subcomplex $\ORGC_d$ spanned by ribbon quivers or oriented ribbon graphs (the ones with no closed paths of directed edges), and it was proven in \cite{Me2} that
$$
H^\bu(\RGC_d)\simeq H^\bu(\ORGC_{d+1}),
$$
again in a close analogy to the commutative case.

\sip

In this note we address the following question: are there low-valence models for
$\ORGC_{d+1}$ (and hence for $H^\bu_c(\cM_{g,m})$) as in the case of $\OGC_{d+1}$? At this point the analogy to the commutative case breaks: such a model generated by graphs with at most {\it tri}valent vertices  is impossible (see \S 2 below and \cite{Me2} where the case $m=1$ was studied). Surprisingly enough,  a low-valence combinatorial model $\sRGC_{d}$ for $H^\bu_c(\cM_{g,m})$ does exist if one relaxes the valency condition for ribbon quivers from {\it at most three}  to {\it at most four}. The precise definition of 
 $\sRGC_{d}$ is given in \S 2 below, and we prove in \S 3 the main theorem of this paper stating the isomorphism
$$
H^\bu(\sRGC_{d+1})=\prod_{g\geq 0, m\geq 1\atop 2g+m\geq 3} H^{\bu +(d-1)m+d(2g-1)}_c(\cM_{g,m}).
$$
of cohomology groups. Thus every cohomology class of the moduli space $\cM_{g,m}$ for any $2g+m\geq 3$, $m\geq 1$ can be represented combinatorially by a  ribbon quiver $\Ga$ with at most four-valent vertices (the {\it at most four}-valency condition is sharp); moreover one can assume without loss of generality that $\Ga$ is invariant under the involution $\iota$ which reverses the direction of its every edge
(see \S {\ref{2: subsec on involution}} for the precise definition of $\iota$). 

  \subsection{Some notation} We work over a field $\K$ of characteristic zero; all the cohomology groups we consider are over $\K$.
 The set $\{1,2, \ldots, n\}$ is abbreviated to $[n]$;  its group of automorphisms is
denoted by $\bS_n$; the trivial (resp., the sign) one-dimensional representation of
 $\bS_n$ is denoted by $\id_n$ (resp.,  $\sgn_n$). The cardinality of a finite set $S$ is 
 denoted by $\# S$ while its linear span over a
field $\K$ by $\K\left\langle S\right\rangle$.
If $V=\oplus_{i\in \Z} V^i$ is a graded vector space, then
$V[k]$ stands for the graded vector space with $V[k]^i:=V^{i+k}$. For 
$v\in V^i$ we set $|v|:=i$. 

\mip

{\bf Acknowledgement}. It is a pleasure to thank Assar Andersson for a useful comment.

\bip

{\Large
\section{\bf Reminder on complexes of ribbon graphs}
}

\mip

\subsection{Penner's ribbon graph complex} Let $\RGC_0$ stand 
for the famous Penner's ribbon graph complex \cite{Pe} which computes the compactly supported cohomology (\ref{1: RGC_0 versus H(M)}) of the moduli spaces $\cM_{g,m}$ of genus $g$ algebraic curves with $m$ marked points.  It can be defined rigorously as follows: consider a set $Rb(m,n,K)$ of connected ribbon graphs with $m$ labelled boundaries, $n$ labelled vertices of valency $\geq 3$ and $k$ labelled undirected edges. The permutation groups $\bS_n\times \bS_k$ acts on this set by relabelling the vertices and the edges. As a $\Z$-graded vector space, $\RGC_0$ is defined as follows,
$$
\RGC_0:= \prod_{m\geq 1,n\geq 1} \K\left\langle Rb(m,n,k)\right\rangle\ot_{\bS_n\times \bS_k} (\id_n \ot \sgn_k)[-k]
$$
Thus a generator $\Ga$ has cohomological degree $\# E(\Ga)$, where $E(\Ga)$ is the set of edges, and some ordering of its edges is assumed to be chosen (up to a permutation and the multiplication by the sign of that permutation); this choice is called an {\it orientation}\ $or$ of $\Ga$. Here are some examples of generators of $\RGC_0$
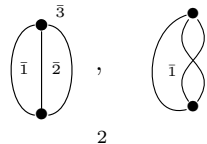
$$
 \Ba{c}\resizebox{12mm}{!}{
\xy
(2.0,-3.5)*{^{{^{\bar{2}}}}},
 (-2,-3.5)*{^{{^{\bar{1}}}}},
(2.5,3.6)*{^{{^{\bar{3}}}}},
 (0,-8)*{\bu}="C";
(0,3)*{\bu}="A1";
(0,3)*{\bu}="A2";
"C"; "A1" **\crv{(-5,-9) & (-5,4)};
"C"; "A2" **\crv{(5,-9) & (5,4)};
 \ar @{-} "A1";"C" <0pt>
\endxy}
\Ea,
%%%%%
 \ \Ba{c}\resizebox{11mm}{!}{ \xy
  (-2.5,-0.5)*{^{{^{\bar{1}}}}},
 (0,8)*{\bu}="1";
%(0,5)*{\circ}="1";
(0,-4)*{\bu}="3";
"1";"3" **\crv{(-5,2) & (5,2)};
"1";"3" **\crv{(5,2) & (-5,2)};
"1";"3" **\crv{(-7,7) & (-7,-7)};
%\ar @{-} "1";"3" <0pt>
\endxy}\Ea 
$$
with three labelled boundaries and one labelled boundary respectively.
The differential in  $\RGC_0$ is given on an arbitrary generator $\Ga\in \RGC_0$ 
by splitting its vertices
\Beq\label{2: delta in RGC_0}
\delta\Ga= (-1)^{|\Ga|}\sum_{v\in V_\bu(\Ga)} \Ga\circ_v  \left(\xy
 (0,0)*{\bu}="a",
(5,0)*{\bu}="b",
\ar @{-} "a";"b" <0pt>
\endxy\right).
\Eeq
where 
$\Ga\circ_v  \left(\xy
 (0,0)*{\bu}="a",
(5,0)*{\bu}="b",
\ar @{-} "a";"b" <0pt>
\endxy\right)$ is a linear combination of ribbon graphs 
obtained from $\Ga$ by substituting into the vertex $v$ the ribbon graph $\xy
 (0,0)*{\bu}="a",
(5,0)*{\bu}="b",
\ar @{-} "a";"b" <0pt>
\endxy$ and then taking a sum over all possible reattachments of the edges (attached earlier to $v$) among the two newly created vertices in a way which respects cyclic order of edges and makes each new  vertex at least trivalent. In particular, if $v$ is itself trivalent, then $\Ga\circ_v  \left(\xy
 (0,0)*{\bu}="a",
(5,0)*{\bu}="b",
\ar @{-} "a";"b" <0pt>
\endxy\right)=0$. 

\sip

This differential preserves the number of boundaries $m$ of a ribbon graph $\Ga\in \RGC_0$ and its genus defined by
$$
 g= 1+\frac{1}{2}\left(\# E(\Ga) - \# V(\Ga)- \# B(\Ga)\right),
 $$
where the symbols  $E(\Ga)$, $V(\Ga)$ and $B(\Ga)$ stand, respectively, for the sets of edges, vertices and boundaries of $\Ga$. Hence Penner's complex decomposes into a direct product\footnote{From now on we work only with ribbon graphs which satisfy the stability condition $2g+m\geq 3$. }
$$
\RGC_0= \prod_{g\geq 0, m\geq 1\atop 2g+m\geq 3}\RGC_0^{g,m}
$$
where $\RGC_0^{g,m}$ is spanned by genus $g$ ribbon graphs with $m$ labelled boundaries. Here are some non-zero examples (in which we skip showing labels of boundaries)
$$
\Ba{c}\mbox{\xy
(0,-2)*{\bu}="A";
(0,-2)*{\bu}="B";
"A"; "B" **\crv{(6,6) & (-6,6)};
"A"; "B" **\crv{(6,-10) & (-6,-10)};
\endxy}\Ea
\in \RGC_0^{0,3}, \ \ \ \ \xy
(0,5)*{\bu}="1";
(0,-4)*{\bu }="3";
%(0,-4)*{\circ}="3";
"1";"3" **\crv{(-5,2) & (5,2)};
"1";"3" **\crv{(5,2) & (-5,2)};
"1";"3" **\crv{(-7,7) & (-7,-7)};
%\ar @{-} "1";"3" <0pt>
\endxy\in \RGC_0^{1,1},\ \ \ \ \Ba{c}\resizebox{9mm}{!}{ \xy
(0,5)*{\bu}="1";
(0,-4)*{\bu}="3";
"1";"3" **\crv{(4,0) & (4,1)};
"1";"3" **\crv{(-4,0) & (-4,-1)};
\ar @{-} "1";"3" <0pt>
\endxy}\Ea \in \RGC_0^{0,3},
$$
while the graph with one boundary
$
\xy
(0,-2)*{\bu}="A";
(0,-2)*{\bu}="B";
"A"; "B" **\crv{(-4,6) & (10,6)};
"A"; "B" **\crv{(-10,6) & (4,6)};
\endxy
$
vanishes identically in $\RGC_0^{1,1}$ as it admits an automorphism which changes
the ordering of its odd edges, i.e.\  it equals minus itself.
One has \cite{Pe, Ko0}
$$
H^\bu(\RGC_0^{g,m})= H^{\bu -m}_c(\cM_{g,m}).
$$

\subsection{Directed (degree shifted) version of Penner'complex}  Consider a  set $Rb^\uparrow(m,n,k)$  of connected ribbon graphs with $m$ labelled boundaries, $n$ labelled vertices of valency  $\geq 2$ (and at least one vertex has valency $\geq 3$) and $k$ labelled {\it directed edges}. The permutation groups $\bS_n\times \bS_k$ acts on this set by relabelling the boundaries and, respectively, the edges. Fix any integer
$d\in \Z$, define a $\Z$-graded vector space,
$$
\dRGC_d:= \prod_{m\geq 1,n\geq 1} \K\left\langle Rb^\uparrow(m,n,k)\right\rangle\ot_{\bS_n\times \bS_k} (\sgn_n^{\ot |d|} \ot \sgn_k^{\ot |d+1|})[k(d-1)],
$$
 An element $\dRGC_d$ can be understood as a pair $(\Ga,or)$ consisting of a ribbon graph $\Ga$ 
whose (directed) edges and vertices are unlabelled, and of a choice $or$ of an ordering of its vertices (for $d$ odd) or an ordering of its edges (for $d$ even) up to an even permutation. There are two possible orientations $or$ and $-or$ on each $\Ga$, and one identifies $(\Ga,-or)=-(\Ga,or)$. Abusing the notation, a pair $(\Ga,or)$ is  abbreviated to $\Ga$ with a choice of $or$ being tacitly assumed.

\sip

The graded vector space $\dRGC_d$ can be made into a complex with the differential
given by a standard ``splitting" formula, 
\Beq\label{5: delta in PreCYd}
\delta\Ga:= \sum_{v\in V(\Ga)} \Ga\circ_v  \left(\xy
 (0,0)*{\bu}="a",
(5,0)*{\bu}="b",
\ar @{->} "a";"b" <0pt>
\endxy\right) \ \ \ \forall \Ga\in \dRGC_d,
\Eeq
in close analogy to (\ref{2: delta in RGC_0}). By the same analogy one has
a decomposition of the complex $\dRGC_d$ into a direct product of subcomplexes
$$
\dRGC_d= \prod_{g\geq 0, m\geq 1\atop 2g+m\geq 3}\dRGC_d^{g,m}.
$$
There is a quasi-isomorphism of complexes
$$
\RGC_0 \lon \dRGC_0
$$
which sends a ribbon graph $\Ga\in \RGC_0$  into a sum of graphs obtained from $\Ga$ by choosing directions on edges in all possible ways (cf.\ Proposition K.1 in \cite{Wi1}).

\sip

For different $d$ the complexes
$\dRGC_d^{g,m}$ are all isomorphic to each other (see e.g.\ \S 3.4 in \cite{Me0}),
$$
 \dRGC_d^{g,m}=\dRGC^{g,m}_0[d(2g-1+m)]
$$
so that one has
$$
H^\bu(\dRGC_d^{g,m})= H^{\bu +(d-1)m+d(2g-1)}_c(\cM_{g,m}).
$$
The cohomological degree of a generator $\Ga\in \dRGC_d$
is given by a formula
$$
|\Ga|=d(\# V(\Ga)-1) + (1-d)\# E(\Ga),
$$
which is identical to the one used for the generators of the Kontsevich graph complex $\GC_d$.

\sip

One can consider a version $\RGC_d$ of $\dRGC_d$ in which the directions of edges are not fixed, but  can be flipped with the following sign factor 
$$
\xy
 (0,0)*{\bu}="a",
(5,0)*{\bu}="b",
\ar @{->} "a";"b" <0pt>
\endxy
=
(-1)^d
\xy
 (0,0)*{\bu}="a",
(5,0)*{\bu}="b",
\ar @{<-} "a";"b" <0pt>
\endxy.
$$
One has $H^\bu(\RGC_d)=H^\bu(\dRGC_d)$ in the full analogy with the case $d=0$.

\sip

We can and shall assume from now on that the generators of the  complex 
$\dRGC_d$ do not containing the so called {\it passing}\, vertices, that is, bivalent vertices with one incoming directed edge and one outgoing directed edge. The directed ribbon graphs having at least one passing vertex span an acyclic subcomplex and hence can be ignored.

\subsection{A complex of ribbon quivers $\ORGC_d$} The complex $\dRGC_d$ contains a subcomplex $\ORGC_d$ spanned by directed ribbon graphs $\Ga$ such that
 the directed edges of $\Ga$ never form closed directed paths, for example
$$
\Ba{c}\resizebox{13mm}{!}{
\xy
   {\ar@/^0.6pc/(-5,0)*{};(5,0)*{\bu}};
 {\ar@{<-}@/^0.6pc/(5,0)*{\bu};(-5,0)*{\bu}};
 (-5,0)*{}="a",
(5,0)*{}="b",
\ar @{->} "a";"b"
\endxy}\Ea\in \ORGC_d , 
\ \ \
\Ba{c}\resizebox{13mm}{!}{
\xy
%(-5,2)*{^1},
%(5,2)*{^2},
%
   {\ar@/^0.6pc/(-5,0)*{};(5,0)*{\bu}};
 {\ar@{<-}@/^0.6pc/(5,0)*{\bu};(-5,0)*{\bu}};
 (-5,0)*{}="a",
(5,0)*{}="b",
\ar @{<-} "a";"b"
\endxy}\Ea\notin \ORGC_d.
$$
It was proven in \cite{Me2} that
$$
H^\bu(\ORGC_{d+1})\simeq  H^\bu(\RGC_d).
$$
This is analogous to a result in the theory of ``commutative"  graph complexes where one has an isomorphism \cite{Wi2, Z, MWW} 
$$
H^\bu(\GC_d)=H^\bu(\OGC_{d+1}).
$$
We refer, e.g., to \S 2 in \cite{MWW} for a short and self-contained definition and description of the ``commutative" complexes $\GC_d$ and $\OGC_d$; they both can be quickly understood via the natural epimorphisms of graded vector spaces,
$$
\RGC_d\rar \GC_d, \ \ \ \ \ORGC_d\rar \OGC_d
$$
which forget cyclic ordering of edges at each vertex of any generator on the left hand side.

\sip

Generators of $\ORGC_{d+1}$ are called {\it oriented}\, ribbon graphs or {\it ribbon quivers}. We use both terms interchangeably. The complex $\ORGC_{d+1}$ decomposes (as usually) into the direct product of subcomplexes
$$
\ORGC_{d+1}= \prod_{g\geq 0, m\geq 1\atop 2g+m\geq 3}\ORGC_{d+1}^{g,m}.
$$
A vertex of a ribbon graph $\Ga$ from $\ORGC_{d+1}$ is called a {\it source}
(resp., {\it target}) if it has no attached incoming (resp., outgoing) edges. Every oriented ribbon graph has at least one source and at least one target.

\subsection{A flow-reversing involution of $\ORGC_d$}\label{2: subsec on involution} Given any ribbon quiver $\Ga\in  \ORGC_d$, there is an associated ribbon quiver $\Ga^{op}\in \ORGC_d$
obtained from $\Ga$ by reversing the direction of its every edges while keeping its orientation $or$ unchanged.
In the full analogy to the commutative case \cite{MZ},  there is an automorphism $\iota$ 
of the complex $\ORGC_d$ given by (cf.\ \cite{MZ})
\Beq\label{2: involution iota}
\Ba{rccl}
\iota: & \ORGC_d & \lon &\ORGC_d\\
&   \Ga & \lon &   \iota(\Ga):=\left\{
\Ba{ll}
(-1)^{\# V(\Ga)+\# E(\Ga) +1}\Ga^{op} & \text{if $d$ is even},\\
(-1)^{\# V(\Ga) +1}\Ga^{op} & \text{if $d$ is odd.}
\Ea
\right. 
\Ea
\Eeq
Hence the complex $\ORGC_d$ decomposes into the direct sum of subcomplexes 
$$
\ORGC_d=\ORGC_d^+ \oplus \ORGC_d^-
$$
corresponding to subcomplexes corresponding to eigenvalues $+1$ and $-1$ of $\iota$ respectively.
There is a similar decomposition in the commutative case 
$$
\OGC_d=\OGC_d^+ \oplus \OGC_d^-
$$
and it was proven in \S 3  of \cite{MZ} that the subcomplex $\OGC_d^-$
is acyclic. That proof is based on a filtration of the complex $\OGC_d$ by the number of vertices of with valencies $\geq 3$ (so that only the part of the differential creating bivalent sources and bivalent targets contributes) and showing that the associated graded $gr \OGC_d^-$ is acyclic. The argument \S 3 of  \cite{MZ})  works {\it word by word}\, for ribbon quivers and leads us to the following conclusion.

\subsubsection{\bf Proposition}\label{2: Prop on aceclicity ORGC_minus} {\it The complex $\ORGC_d^-$ is acyclic so that
$H^\bu(\ORGC_d)=H^\bu(\ORGC^+_d)$.}

\subsection{Remark on labellings of boundaries}\label{ 2: subsec on S_m}   The permutation ghroup $\bS_m$ acts on the complexes $\RGC_d^{g,m}$ and $\ORGC_d^{g,m}$ by relabelling the boundaries of ribbon graphs. One must be a bit careful when translating this action to the action on the marked points of the associated cohomology classes of $H^\bu(\cM_{g,m})$. Indeed, the isomorphism (\ref{1: RGC_0 versus H(M)}) says that we have the following {\it isomorphism of $\bS_m$-modules},
\Beq\label{3: RGC_0 versus H(M) as S_m modules}
H^\bu(\RGC_0^{g,m})= H^{\bu}_c(\cM_{g,m} \times (\R^{>0})^m)= H^{\bu -m}_c(\cM_{g,m}) \ot\sgn_m
\Eeq                                                        
where in the right most term a diagonal action of $\bS_m$ on the tensor product is assumed, and $\R^+$ stands for the space of positive real numbers. This implies, 
for example, that
$$
 H^3(\RGC_0^{0,3})= H^0(\cM_{0,3})\ot \sgn_3\simeq \sgn_3,
$$
Hence the cohomology class of the point $\cM_{0,3}$ is expected to be represented in $\RGC_0^{0,3}$  by a cycle  with {\it 
 skew-symmetrized}\, labelled boundaries.
 
 \sip
 
It is worth emphasizing in this context that the above mentioned isomorphism of cohomology groups
 $$
   H^\bu(\ORGC_{d+1}^{g,m})=\RGC_{d}^{g,m}
 $$
is in fact an {\it isomorphism of $\bS_m$-modules}. In particular, one must have an equality
    \Beq
 H^3(\ORGC_1^{0,3})= H^0(\cM_{0,3})\ot \sgn_3\simeq sgn_3.
 \Eeq     
This may seem surprising as the orientation of ribbon graphs are defined 
differently for $d$ even and $d$ odd (see \S 3.4 in \cite{Me0} 
for more details on the role of the integer parameter $d$).

\subsection{On complexes of oriented graphs with at most trivalent vertices} Let $I$ be a subspace 
in the commutative graph complex $\OGC_d$ spanned by graphs having at least one vertex of valency $\geq 4$ or at least one trivalent source or trivalent target, and let $\langle I, \delta I\rangle$ be its differential closure. There is a short exact sequence of complexes
$$
0 \lon \langle I, \delta I\rangle \lon \OGC_{d+1} \stackrel{p}{\lon} \OGC_{d+1}^T \lon 0
$$
where the quotient complex $\OGC_{d+1}^T$ is generated by equivalence classes of graphs whose vertices can be only of the form 
\Beq\label{2: bivalent abd trivalent vertices}
\Ba{c}\resizebox{8mm}{!}{ \xy
 (0,0)*{}="a",
(4,3)*{\bullet}="b",
(8,0)*{}="c",
\ar @{<-} "a";"b" <0pt>
\ar @{->} "b";"c" <0pt>
\endxy}\Ea, 
\ \ \ \  
\Ba{c}\resizebox{8mm}{!}{\xy
 (0,0)*{}="a",
(4,3)*{\bullet}="b",
(8,0)*{}="c",
\ar @{->} "a";"b" <0pt>
\ar @{<-} "b";"c" <0pt>
\endxy} \Ea, 
\ \ \ 
 \Ba{c}\resizebox{8mm}{!}{  \xy
(0,6)*{}="1";
    (0,0.2)*{\bu}="L";
  (-4,-5)*{}="C";
   (+4,-5)*{}="D";
\ar @{->} "D";"L" <0pt>
\ar @{->} "C";"L" <0pt>
\ar @{<-} "1";"L" <0pt>
 \endxy}
 \Ea,
 \ \ \
 \Ba{c}\resizebox{8mm}{!}{  \xy
(0,-6)*{}="1";
    (0,-0.2)*{\bu}="L";
  (-4,5)*{}="C";
   (4,5)*{}="D";
\ar @{<-} "D";"L" <0pt>
\ar @{<-} "C";"L" <0pt>
\ar @{->} "1";"L" <0pt>
 \endxy}
 \Ea
\Eeq

They are subject to (co)Jacobi or Drinfeld type equivalence relations for every edge connecting two such vertices, and the induced differential on   $\OGC_{d+1}^T$ acts only on trivalent vertices by creating new bivalent ones. It was proven in \cite{Me1} that the epimorphism
$$
p: \OGC_{d+1}\rar \OGC_{d+1}^T
$$
is a quasi-isomorphism, so that one has
$$
H^\bu(\GC_d)\simeq H^\bu(\OGC_{d+1}^T).
$$
Thus every cohomology class in the Kontsevich graph complex $\GC_d$ can be represented in terms of quivers with 2 and 3 valent vertices only.

\sip

One can repeat the above scenario for the oriented ribbon graph complex, 
$$
0 \lon \langle I, \delta I\rangle \lon \ORGC_{d+1} \stackrel{\pi}{\lon} \ORGC_{d+1}^T \lon 0,
$$
and obtain a quotient complex $\ORGC_{d+1}^T$ spanned by equivalence class 
of ribbon quivers whose vertices are of the form (\ref{2: bivalent abd trivalent vertices}) only. However in this case the epimorphism 
$$
\ORGC_{d+1} \stackrel{\pi}{\lon} \ORGC_{d+1}^T 
$$
is {\it not}\, a quasi-isomorphism, and hence it is {\it not}\, true that 
 every cohomology class in $H^\bu_c(\cM_{g,m})$ can be represented in terms of quivers with 2 and 3 valent vertices only. Indeed, if a cohomology class $[\Ga]\in H^k_c(\cM_{g,m})$ of degree $k$ is represented by a cycle $\Ga \in \ORGC_{d+1}$ with $p_2$ bivalent vertices and $p_3$ trivalent vertices of type (\ref{2: bivalent abd trivalent vertices}), then one must the have equality
 $$
 k= |\Ga| + 2dg - m + dm - d =p_2+ 4g-5 +m  
 $$
 which implies $k\geq 4g-3+m$ as $p_2\geq 2$.  There are many well-known non-vanishing  results for cohomology groups  $H^k_c(\cM_{g,m})$ with  $k< 4g-3+m$, e.g.\  $H^4_c(\cM_{2,1}) \neq 0$. Hence it can {\it not}\,  be true that 
  every cohomology class of the moduli space $\cM_{g,m}$ for any $2g+m\geq 3$, $m\geq 1$ can be represented combinatorially by a ribbon quiver with at most three-valent vertices. 

\sip

 Surprisingly, it all works fine if we add to the  list (\ref{2: bivalent abd trivalent vertices}) of allowed vertices the  4-valent vertex of the following  special form:
 \Beq\label{2: fourvalent vertex}
  \Ba{c}\resizebox{13mm}{!}{  \xy
(0,-5)*{}="U";
(0,+5)*{}="D";
    (0,0)*{\bu}="C";
  (5,0)*{}="R";
   (-5,0)*{}="L";
\ar @{->} "C";"U" <0pt>
\ar @{->} "C";"D" <0pt>
\ar @{<-} "C";"L" <0pt>
\ar @{<-} "C";"R" <0pt>
 \endxy}
 \Ea
\Eeq
 This is the main point of our paper.

\subsection{A complex of ribbon quivers with at most four-valent vertices}
Let $J$ be a linear subspace 
of the ribbon graph complex $\ORGC_{d+1}$ spanned by graphs having at least one trivalent source or trivalent target, or at least one vertex of valency $\geq 4$ which is {\it not}\, of the form (\ref{2: fourvalent vertex}), and let $\langle J, \delta J\rangle$ be its differential closure. Define a low-valence ribbon graph complex by a short exact sequence of complexes
$$
0 \lon \langle J, \delta J\rangle \lon \ORGC_{d+1} \stackrel{p}{\lon} \sRGC_{d+1} \lon 0.
$$
The quotient complex  $\sRGC_{d+1}$ is generated by equivalence classes of graphs whose vertices have valency at most four\footnote{The symbol $\Delta$ stands for number 4 in Greek.} and can be only of one of the following four types 
\Beq\label{2: types of vertices in sORGC}
\Ba{c}\resizebox{8mm}{!}{ \xy
 (0,0)*{}="a",
(4,3)*{\bullet}="b",
(8,0)*{}="c",
\ar @{<-} "a";"b" <0pt>
\ar @{->} "b";"c" <0pt>
\endxy}\Ea, 
\ \ \ \  
\Ba{c}\resizebox{8mm}{!}{\xy
 (0,0)*{}="a",
(4,3)*{\bullet}="b",
(8,0)*{}="c",
\ar @{->} "a";"b" <0pt>
\ar @{<-} "b";"c" <0pt>
\endxy} \Ea, 
\ \ \ 
 \Ba{c}\resizebox{8mm}{!}{  \xy
(0,6)*{}="1";
    (0,0.2)*{\bu}="L";
  (-4,-5)*{}="C";
   (+4,-5)*{}="D";
\ar @{->} "D";"L" <0pt>
\ar @{->} "C";"L" <0pt>
\ar @{<-} "1";"L" <0pt>
 \endxy}
 \Ea,
 \ \ \
 \Ba{c}\resizebox{8mm}{!}{  \xy
(0,-6)*{}="1";
    (0,-0.2)*{\bu}="L";
  (-4,5)*{}="C";
   (4,5)*{}="D";
\ar @{<-} "D";"L" <0pt>
\ar @{<-} "C";"L" <0pt>
\ar @{->} "1";"L" <0pt>
 \endxy}
 \Ea, 
 \ \ \  
  \Ba{c}\resizebox{13mm}{!}{  \xy
(0,-5)*{}="U";
(0,+5)*{}="D";
    (0,0)*{\bu}="C";
  (5,0)*{}="R";
   (-5,0)*{}="L";
\ar @{->} "C";"U" <0pt>
\ar @{->} "C";"D" <0pt>
\ar @{<-} "C";"L" <0pt>
\ar @{<-} "C";"R" <0pt>
 \endxy}
 \Ea\ .
\Eeq
For every edge connecting vertices of the above types there is an equivalence relation (a kind of non-commutative IHX relation) given in terms of {\it ribbon}\,  graphs with labelled legs (standing for edges connected to some other vertex of vertices) as follows:
\sip

{\bf (I)}``2+3"  relations,

\Beq\label{2: 2+3 relations}
 \Ba{c}\resizebox{13.5mm}{!}{  \xy
    (-9,+3)*{\bbu}="L";
 (-14,-3.5)*{\bbu}="B";
 (-20,-12)*+{_1}="b1";
 (-8,-12)*+{_2}="b2";
  (-3,-4)*{_3}="C";
\ar @{->} "C";"L" <0pt>
\ar @{->} "B";"L" <0pt>
\ar @{<-} "B";"b1" <0pt>
\ar @{<-} "B";"b2" <0pt>
 \endxy}
 \Ea
 +
  \Ba{c}\resizebox{13.5mm}{!}{  \xy
    (-9,+3)*{\bbu}="L";
 (-14,-3.5)*{\bbu}="B";
 (-20,-12)*+{_2}="b1";
 (-8,-12)*+{_3}="b2";
  (-3,-4)*{_1}="C";
\ar @{->} "C";"L" <0pt>
\ar @{->} "B";"L" <0pt>
\ar @{<-} "B";"b1" <0pt>
\ar @{<-} "B";"b2" <0pt>
 \endxy}
 \Ea
 +
  \Ba{c}\resizebox{13.5mm}{!}{  \xy
    (-9,+3)*{\bbu}="L";
 (-14,-3.5)*{\bbu}="B";
 (-20,-12)*+{_3}="b1";
 (-8,-12)*+{_1}="b2";
  (-3,-4)*{_2}="C";
\ar @{->} "C";"L" <0pt>
\ar @{->} "B";"L" <0pt>
\ar @{<-} "B";"b1" <0pt>
\ar @{<-} "B";"b2" <0pt>
 \endxy}
 \Ea
 =0\ , 
 %%%%%%%%%%%%
 \hspace{-1mm}
\ \ \ 
  \Ba{c}\resizebox{13.5mm}{!}{  \xy
    (-9,+3)*{\bbu}="L";
 (-14,-3.5)*{\bbu}="B";
 (-20,-12)*+{_1}="b1";
 (-8,-12)*+{_2}="b2";
  (-3,-4)*{_3}="C";
\ar @{<-} "C";"L" <0pt>
\ar @{<-} "B";"L" <0pt>
\ar @{->} "B";"b1" <0pt>
\ar @{->} "B";"b2" <0pt>
 \endxy}
 \Ea
 +
  \Ba{c}\resizebox{13.5mm}{!}{  \xy
    (-9,+3)*{\bbu}="L";
 (-14,-3.5)*{\bbu}="B";
 (-20,-12)*+{_2}="b1";
 (-8,-12)*+{_3}="b2";
  (-3,-4)*{_1}="C";
\ar @{<-} "C";"L" <0pt>
\ar @{<-} "B";"L" <0pt>
\ar @{->} "B";"b1" <0pt>
\ar @{->} "B";"b2" <0pt>
 \endxy}
 \Ea
 +
   \Ba{c}\resizebox{13.5mm}{!}{  \xy
    (-9,+3)*{\bbu}="L";
 (-14,-3.5)*{\bbu}="B";
 (-20,-12)*+{_3}="b1";
 (-8,-12)*+{_1}="b2";
  (-3,-4)*{_2}="C";
\ar @{<-} "C";"L" <0pt>
\ar @{<-} "B";"L" <0pt>
\ar @{->} "B";"b1" <0pt>
\ar @{->} "B";"b2" <0pt>
 \endxy}
 \Ea
 \hspace{-1mm}
 =0\ , 
 \hspace{-1mm}
\Eeq
%%%%%%%%%%%%%%%%%%%%%%%%%%%%%%%%%%%%%%%%%%%%%%%%
{\bf (II)} twisted {cyclic (co)associativity relations},
\Beq\label{2: twisted ass-type relations}
 \Ba{c}\resizebox{13.5mm}{!}{  \xy
(-9,10)*{_0}="1";
    (-9,+3)*{\bbu}="L";
 (-14,-3.5)*{\bbu}="B";
 (-20,-12)*+{_1}="b1";
 (-8,-12)*+{_2}="b2";
  (-3,-4)*{_3}="C";
\ar @{<-} "C";"L" <0pt>
\ar @{<-} "B";"L" <0pt>
\ar @{->} "B";"b1" <0pt>
\ar @{->} "B";"b2" <0pt>
\ar @{->} "1";"L" <0pt>
 \endxy}
 \Ea
 +
  \Ba{c}\resizebox{13.5mm}{!}{  \xy
(9,10)*{_0}="1";
    (9,+3)*{\bbu}="L";
 (14,-3.5)*{\bbu}="B";
 (20,-12)*+{_3}="b1";
 (8,-12)*+{_2}="b2";
  (3,-4)*{_1}="C";
\ar @{<-} "C";"L" <0pt>
\ar @{<-} "B";"L" <0pt>
\ar @{->} "B";"b1" <0pt>
\ar @{->} "B";"b2" <0pt>
\ar @{->} "1";"L" <0pt>
 \endxy}
 \Ea
 -
    \Ba{c}\resizebox{15mm}{!}{  \xy
 (0,-14)*{_2}="1";
(0,7)*{^0}="U";
(0,-7)*{\bu}="D";
    (0,0)*{\bu}="C";
  (7,0)*{_3}="R";
   (-7,0)*{_1}="L";
\ar @{<-} "C";"U" <0pt>
\ar @{<-} "C";"D" <0pt>
\ar @{->} "C";"L" <0pt>
\ar @{->} "C";"R" <0pt>
\ar @{->} "D";"1" <0pt>
 \endxy}
 \Ea
 =0, \ \ 
 \Ba{c}\resizebox{13.5mm}{!}{  \xy
(-9,10)*{^0}="1";
    (-9,+3)*{\bbu}="L";
 (-14,-3.5)*{\bbu}="B";
 (-20,-12)*+{_1}="b1";
 (-8,-12)*+{_2}="b2";
  (-3,-4)*{_3}="C";
\ar @{->} "C";"L" <0pt>
\ar @{->} "B";"L" <0pt>
\ar @{<-} "B";"b1" <0pt>
\ar @{<-} "B";"b2" <0pt>
\ar @{<-} "1";"L" <0pt>
 \endxy}
 \Ea
 +
  \Ba{c}\resizebox{13.5mm}{!}{  \xy
(9,10)*{^0}="1";
    (9,+3)*{\bbu}="L";
 (14,-3.5)*{\bbu}="B";
 (20,-12)*+{_3}="b1";
 (8,-12)*+{2}="b2";
  (3,-4)*{1}="C";
\ar @{->} "C";"L" <0pt>
\ar @{->} "B";"L" <0pt>
\ar @{<-} "B";"b1" <0pt>
\ar @{<-} "B";"b2" <0pt>
\ar @{<-} "1";"L" <0pt>
 \endxy}
 \Ea
 -
   \Ba{c}\resizebox{15mm}{!}{  \xy
 (0,-14)*{_2}="1";
(0,7)*{^0}="U";
(0,-7)*{\bu}="D";
    (0,0)*{\bu}="C";
  (7,0)*{_3}="R";
   (-7,0)*{_1}="L";
\ar @{->} "C";"U" <0pt>
\ar @{->} "C";"D" <0pt>
\ar @{<-} "C";"L" <0pt>
\ar @{<-} "C";"R" <0pt>
\ar @{<-} "D";"1" <0pt>
 \endxy}
 \Ea
 \hspace{-1mm} =0, 
\Eeq

\sip

{\bf (III)} {\sf infinitesimal bialgebra relation},
\Beq\label{2: IB type relation}
  \Ba{c}\resizebox{6.5mm}{!}{  \xy
(-4,10)*{}="1";
 (4,10)*{}="2";
    (0,3.5)*{\bbu}="A";
 (0,-3.5)*{\bbu}="B";
 (-4,-10)*{}="b1";
 (4,-10)*{}="b2";
\ar @{->} "A";"1" <0pt>
\ar @{->} "A";"2" <0pt>
\ar @{->} "B";"A" <0pt>
\ar @{<-} "B";"b1" <0pt>
\ar @{<-} "B";"b2" <0pt>
 \endxy}
 \Ea
 +
 \Ba{c}\resizebox{13.5mm}{!}{  \xy
(-9,10)*{}="1";
    (-9,+3)*{\bbu}="L";
 (-14,-3.5)*{\bbu}="B";
 (-19,5)*+{}="b1";
 (-14,-12)*+{}="b2";
  (-3,-5)*{}="C";
\ar @{->} "C";"L" <0pt>
\ar @{->} "B";"L" <0pt>
\ar @{->} "B";"b1" <0pt>
\ar @{<-} "B";"b2" <0pt>
\ar @{<-} "1";"L" <0pt>
 \endxy}
 \Ea
 +
 \Ba{c}\resizebox{13.5mm}{!}{  \xy
(-9,-10)*{}="1";
    (-9,-3)*{\bbu}="L";
 (-14,3.5)*{\bbu}="B";
 (-19,-5)*+{}="b1";
 (-14,12)*+{}="b2";
  (-3,5)*{}="C";
\ar @{<-} "C";"L" <0pt>
\ar @{<-} "B";"L" <0pt>
\ar @{<-} "B";"b1" <0pt>
\ar @{->} "B";"b2" <0pt>
\ar @{->} "1";"L" <0pt>
 \endxy}
 \Ea
=0
\Eeq

{\bf (IV)} {``$3+4$"-relations},
\Beq\label{2: 3+4 relations}
  \Ba{c}\resizebox{17mm}{!}{  \xy
(-3,13)*+{_1}="UL";
(3,13)*+{_2}="UR";
(0,7)*{\bu}="U";
(0,-7)*+{_4}="D";
    (0,0)*{\bu}="C";
  (7,0)*+{_3}="R";
   (-7,0)*+{_0}="L";
\ar @{->} "U";"UL" <0pt>   
\ar @{->} "U";"UR" <0pt>  
\ar @{->} "C";"U" <0pt>
\ar @{->} "C";"D" <0pt>
\ar @{<-} "C";"L" <0pt>
\ar @{<-} "C";"R" <0pt>
 \endxy}
 \Ea
 \hspace{-1mm}
 -
 \hspace{-2mm}
   \Ba{c}\resizebox{17mm}{!}{  \xy
(-3,13)*+{_0}="UL";
(3,13)*+{_1}="UR";
(0,7)*{\bu}="U";
(0,-7)*+{_3}="D";
    (0,0)*{\bu}="C";
  (7,0)*+{_2}="R";
   (-7,0)*+{_4}="L";
\ar @{<-} "U";"UL" <0pt>   
\ar @{->} "U";"UR" <0pt>  
\ar @{<-} "C";"U" <0pt>
\ar @{<-} "C";"D" <0pt>
\ar @{->} "C";"L" <0pt>
\ar @{->} "C";"R" <0pt>
 \endxy}
 \Ea
 \hspace{-1mm}
 -
 \hspace{-2mm}
    \Ba{c}\resizebox{17mm}{!}{  \xy
(-3,13)*+{_2}="UL";
(3,13)*+{_3}="UR";
(0,7)*{\bu}="U";
(0,-7)*+{_0}="D";
    (0,0)*{\bu}="C";
  (7,0)*+{_4}="R";
   (-7,0)*+{_1}="L";
\ar @{->} "U";"UL" <0pt>   
\ar @{<-} "U";"UR" <0pt>  
\ar @{<-} "C";"U" <0pt>
\ar @{<-} "C";"D" <0pt>
\ar @{->} "C";"L" <0pt>
\ar @{->} "C";"R" <0pt>
 \endxy}
 \Ea
 \hspace{-2mm}
 =0, 
 %%%%%%%%%%%%%%%%%
   \Ba{c}\resizebox{17mm}{!}{  \xy
(-3,13)*+{_1}="UL";
(3,13)*+{_2}="UR";
(0,7)*{\bu}="U";
(0,-7)*+{_4}="D";
    (0,0)*{\bu}="C";
  (7,0)*+{_3}="R";
   (-7,0)*+{_0}="L";
\ar @{<-} "U";"UL" <0pt>   
\ar @{<-} "U";"UR" <0pt>  
\ar @{<-} "C";"U" <0pt>
\ar @{<-} "C";"D" <0pt>
\ar @{->} "C";"L" <0pt>
\ar @{->} "C";"R" <0pt>
 \endxy}
 \Ea
 \hspace{-1mm}
 -
 \hspace{-2mm}
   \Ba{c}\resizebox{17mm}{!}{  \xy
(-3,13)*+{_0}="UL";
(3,13)*+{_1}="UR";
(0,7)*{\bu}="U";
(0,-7)*+{_3}="D";
    (0,0)*{\bu}="C";
  (7,0)*+{_2}="R";
   (-7,0)*+{_4}="L";
\ar @{->} "U";"UL" <0pt>   
\ar @{<-} "U";"UR" <0pt>  
\ar @{->} "C";"U" <0pt>
\ar @{->} "C";"D" <0pt>
\ar @{<-} "C";"L" <0pt>
\ar @{<-} "C";"R" <0pt>
 \endxy}
 \hspace{0mm}
-
\hspace{-2mm}
     \Ba{c}\resizebox{17mm}{!}{  \xy
(-3,13)*+{_2}="UL";
(3,13)*+{_3}="UR";
(0,7)*{\bu}="U";
(0,-7)*+{_0}="D";
    (0,0)*{\bu}="C";
  (7,0)*+{_4}="R";
   (-7,0)*+{_1}="L";
\ar @{<-} "U";"UL" <0pt>   
\ar @{->} "U";"UR" <0pt>  
\ar @{->} "C";"U" <0pt>
\ar @{->} "C";"D" <0pt>
\ar @{<-} "C";"L" <0pt>
\ar @{<-} "C";"R" <0pt>
 \endxy}
 \Ea
 \Ea
 \hspace{-2mm}
 =0,
\Eeq

{\bf (V)} { double Lie relation},
\Beq\label{2: double Lie relations}
  \Ba{c}\resizebox{28mm}{!}{  \xy
(0,-7)*+{_5}="D";
(0,+7)*+{_1}="U";
 (-7,0)*+{_0}="L";
(0,0)*{\bu}="C";
(7,0)*{\bu}="CR";
  (7,-7)*+{_4}="RD";
(7,+7)*+{_2}="RU";
 (14,0)*+{_3}="RR";
\ar @{->} "C";"U" <0pt>
\ar @{->} "C";"D" <0pt>
\ar @{<-} "C";"L" <0pt>
\ar @{<-} "C";"CR" <0pt>
\ar @{<-} "CR";"RU" <0pt>
\ar @{<-} "CR";"RD" <0pt>
\ar @{->} "CR";"RR" <0pt>
 \endxy}
 \Ea
 +
  \Ba{c}\resizebox{28mm}{!}{  \xy
(0,-7)*+{_3}="D";
(0,+7)*+{_5}="U";
 (-7,0)*+{_4}="L";
(0,0)*{\bu}="C";
(7,0)*{\bu}="CR";
  (7,-7)*+{_2}="RD";
(7,+7)*+{_0}="RU";
 (14,0)*+{_1}="RR";
\ar @{->} "C";"U" <0pt>
\ar @{->} "C";"D" <0pt>
\ar @{<-} "C";"L" <0pt>
\ar @{<-} "C";"CR" <0pt>
\ar @{<-} "CR";"RU" <0pt>
\ar @{<-} "CR";"RD" <0pt>
\ar @{->} "CR";"RR" <0pt>
 \endxy}
 \Ea
  +
  \Ba{c}\resizebox{28mm}{!}{  \xy
(0,-7)*+{_1}="D";
(0,+7)*+{_3}="U";
 (-7,0)*+{_2}="L";
(0,0)*{\bu}="C";
(7,0)*{\bu}="CR";
  (7,-7)*+{_0}="RD";
(7,+7)*+{_4}="RU";
 (14,0)*+{_5}="RR";
\ar @{->} "C";"U" <0pt>
\ar @{->} "C";"D" <0pt>
\ar @{<-} "C";"L" <0pt>
\ar @{<-} "C";"CR" <0pt>
\ar @{<-} "CR";"RU" <0pt>
\ar @{<-} "CR";"RD" <0pt>
\ar @{->} "CR";"RR" <0pt>
 \endxy}
 \Ea=0.
\Eeq
The sign rule is especially simple for $d$  odd
(when the edges are even) which is assumed in the above formulae: all the vertices in the relations {\bf (I)-(IV)} are ordered  from the top to the bottom, and in the last relation {\bf (V)} --- from the left to the right.
 It is tacitly assumed  that the implementation of any the above relations in a ribbon quiver $\Ga$ does not create {\it closed}\, paths 
of directed edges in $\Ga$, i.e.\ ribbon  quivers stay ribbon quivers.

\sip

The differential in the complex $\sRGC_{d+1}$ is given by the sum
over its actions on the vertices,
$$
\delta\Ga:= \sum_{v\in V(\Ga)} 
\delta_v\Ga
$$
and the operation $\delta_v$ changes the vertex $v$ of $\Ga$ in accordance 
with its type as follows,
\Beq\label{2 d on 2val}
\delta_v \left(\Ba{c}\resizebox{9mm}{!}{ \xy
 (0,0)*{}="a",
(4,3)*{\bullet}="b",
(8,0)*{}="c",
\ar @{<-} "a";"b" <0pt>
\ar @{->} "b";"c" <0pt>
\endxy}\Ea\right)=0\ ,
 \ \ \ 
\delta_v 
\left(\Ba{c}\resizebox{9mm}{!}{ \xy
 (0,0)*{}="a",
(4,3)*{\bullet}="b",
(8,0)*{}="c",
%>\
\ar @{->} "a";"b" <0pt>
\ar @{<-} "b";"c" <0pt>
\endxy}\Ea 
\right)=0,
\Eeq
%%%%%%%%%%%%%%%%%%%%%%%%%%%%%%
\Beq\label{2: d on 3 val}
\delta_v 
\left( \Ba{c}\resizebox{8mm}{!}{  \xy
(0,6)*{}="1";
    (0,0.2)*{\bu}="L";
  (-4,-5)*{}="C";
   (+4,-5)*{}="D";
\ar @{->} "D";"L" <0pt>
\ar @{->} "C";"L" <0pt>
\ar @{<-} "1";"L" <0pt>
 \endxy}
 \Ea\right)
 =
  \Ba{c}\resizebox{11.5mm}{!}{  \xy
(0,-6)*{}="1";
(9,-0.2)*{}="2";
    (0,-0.2)*{\bu}="L";
  (-4,5)*{}="C";
   (4,5)*{\bu}="D";
\ar @{<-} "D";"L" <0pt>
\ar @{<-} "C";"L" <0pt>
\ar @{->} "1";"L" <0pt>
\ar @{->} "2";"D" <0pt>
 \endxy}
 \Ea
 +
   \Ba{c}\resizebox{11.5mm}{!}{  \xy
(0,-6)*{}="1";
(-9,-0.2)*{}="2";
    (0,-0.2)*{\bu}="L";
  (4,5)*{}="C";
   (-4,5)*{\bu}="D";
\ar @{<-} "D";"L" <0pt>
\ar @{<-} "C";"L" <0pt>
\ar @{->} "1";"L" <0pt>
\ar @{->} "2";"D" <0pt>
 \endxy}
 \Ea, \ \
\delta_v 
\left( \Ba{c}\resizebox{8mm}{!}{  \xy
(0,-6)*{}="1";
    (0,-0.2)*{\bu}="L";
  (-4,5)*{}="C";
   (+4,5)*{}="D";
\ar @{<-} "D";"L" <0pt>
\ar @{<-} "C";"L" <0pt>
\ar @{->} "1";"L" <0pt>
 \endxy}
 \Ea\right)=
  \Ba{c}\resizebox{11.5mm}{!}{  \xy
(0,6)*{}="1";
(8,0.2)*{}="2";
    (0,0.2)*{\bu}="L";
  (-4,-5)*{}="C";
   (4,-5)*{\bu}="D";
\ar @{->} "D";"L" <0pt>
\ar @{->} "C";"L" <0pt>
\ar @{<-} "1";"L" <0pt>
\ar @{<-} "2";"D" <0pt>
 \endxy}
 \Ea
 +
   \Ba{c}\resizebox{11.5mm}{!}{  \xy
(0,6)*{}="1";
(-8,0.2)*{}="2";
    (0,0.2)*{\bu}="L";
  (4,-5)*{}="C";
   (-4,-5)*{\bu}="D";
\ar @{->} "D";"L" <0pt>
\ar @{->} "C";"L" <0pt>
\ar @{<-} "1";"L" <0pt>
\ar @{<-} "2";"D" <0pt>
 \endxy}
 \Ea
\Eeq
%%%%%%%%%%%%%%%%%%%%
\Beq\label{2: d on 4val}
\delta_v\left(  \Ba{c}\resizebox{17.5mm}{!}{  \xy
(0,-7)*+{_3}="D";
(0,+7)*+{_1}="U";
    (0,0)*{\bu}="C";
  (7,0)*+{_2}="R";
   (-7,0)*+{_0}="L";
\ar @{->} "C";"U" <0pt>
\ar @{->} "C";"D" <0pt>
\ar @{<-} "C";"L" <0pt>
\ar @{<-} "C";"R" <0pt>
 \endxy}
 \Ea\right)
 =
  \Ba{c}\resizebox{18.5mm}{!}{  \xy
(0,-7)*+{_3}="D";
(0,5)*{\bu}="U";
(0,+12)*+{_1}="0";
    (0,0)*{\bu}="C";
  (7,0)*+{_2}="R";
   (-7,5)*+{_0}="L";
\ar @{->} "U";"0" <0pt>  
\ar @{->} "C";"U" <0pt>
\ar @{->} "C";"D" <0pt>
\ar @{<-} "U";"L" <0pt>
\ar @{<-} "C";"R" <0pt>
 \endxy}
 \Ea 
 +
  \Ba{c}\resizebox{18.5mm}{!}{  \xy
(0,-7)*+{_1}="D";
(0,5)*{\bu}="U";
(0,+12)*+{_3}="0";
    (0,0)*{\bu}="C";
  (7,0)*+{_0}="R";
   (-7,5)*+{_2}="L";
\ar @{->} "U";"0" <0pt>  
\ar @{->} "C";"U" <0pt>
\ar @{->} "C";"D" <0pt>
\ar @{<-} "U";"L" <0pt>
\ar @{<-} "C";"R" <0pt>
 \endxy}
 \Ea
  +
  \Ba{c}\resizebox{18.5mm}{!}{  \xy
(0,-7)*+{_3}="D";
(0,5)*{\bu}="U";
(0,+12)*+{_1}="0";
    (0,0)*{\bu}="C";
  (7,5)*+{_2}="R";
   (-7,0)*+{_0}="L";
\ar @{->} "U";"0" <0pt>  
\ar @{->} "C";"U" <0pt>
\ar @{->} "C";"D" <0pt>
\ar @{<-} "U";"R" <0pt>
\ar @{<-} "C";"L" <0pt>
 \endxy}
 \Ea 
   +
  \Ba{c}\resizebox{18.5mm}{!}{  \xy
(0,-7)*+{_1}="D";
(0,5)*{\bu}="U";
(0,+12)*+{_3}="0";
    (0,0)*{\bu}="C";
  (7,5)*+{_0}="R";
   (-7,0)*+{_2}="L";
\ar @{->} "U";"0" <0pt>  
\ar @{->} "C";"U" <0pt>
\ar @{->} "C";"D" <0pt>
\ar @{<-} "U";"R" <0pt>
\ar @{<-} "C";"L" <0pt>
 \endxy}
 \Ea 
\Eeq

Assuming that $d$ is even (the choice of a particular value for $d$  plays no role, it is a matter of convenience), the vertices in the right hand sides of the above formulae are ordered from   the top to the bottom.

\subsection{Main Theorem}\label{2: Main theorem} {\it The epimorphism of complexes
$$
 \ORGC_{d+1} \stackrel{\pi}{\lon} \sRGC_{d+1} 
$$
is a quasi-isomorphism.}

\mip

We prove this theorem in the next section.

\mip

The ``small"  complex inherits from $\ORGC_{d+1}$ the ``genus+boundary"  decomposition 
$$
\sRGC_{d+1}= \prod_{g\geq 0, m\geq 1\atop 2g+m\geq 3}\sRGC_{d+1}^{g,m}
$$
as well as the eigenvalue decomposition with respect to the involution (\ref{2: involution iota}), 
$$
 \sRGC_{d+1}= \sRGC_{d+1}^+ \oplus \sRGC_{d+1}^-, \ \ \ \sRGC_{d+1}^{g,m}=
 \sRGC_{d+1}^{+,g,m}\oplus \sRGC_{d+1}^{-,g,m}.
$$
By Proposition {\ref{2: Prop on aceclicity ORGC_minus}}, the subcomplex
 $\sRGC_{d+1}^-$ is acyclic. This observation together with the Main Theorem
imply the following result.

\subsection{Corollary}\label{2: iso of sRGC+,g,m with M} {\it One has isomorphisms of cohomology groups}
\Beq\label{2: Main equality}
H^\bu(\sRGC_{d+1}^{g,m})\simeq H^\bu(\sRGC_{d+1}^{+,g,m})\simeq H^{\bu +(d-1)m+d(2g-1)}_c(\cM_{g,m})\simeq  H^{(3-d)(2g+m) +d-6 -\bu}(\cM_{g,m}).
\Eeq

\subsection{Corollary}\label{2: Corollary on ORCG to dRGC} {\it The inclusion of complexes
\Beq\label{2: ORGC to dRGC}
i: \ORGC_{d} \lon \dRGC_d 
\Eeq
induces the zero map at the cohomology level for any $d\in \Z$.}

\begin{proof}
Any cycle representing of a cohomology class in $H^\bu(\ORGC_d)$ must contain summands with at least one bivalent source and at least one bivalent target. There are no such cohomology classes in $H^\bu(\dRGC_d)$ as, thanks to the quasi-isomorphism (\ref{1: RGC to dRGC}), they are given by ribbon graphs with every vertex at least trivalent.
\end{proof}

A similar result holds true \cite{Z2} in the commutative case due to (i) the lower bounds  on degrees of non-vanishing classes stemming from the at least trivalency conditions on classes from $H^\bu(\dGC_d)$, and (ii)  the upper bounds computed in \cite{Wi4}. The vanishing of the morphism $H^\bu(\OGC_d) \rar H^\bu(\dGC_d)$ also follows from Proposition 3.2.1 in \cite{Me2}.

\subsection{Compatibility with Harer's vanishing theorem} The only known at present ``universal" vanishing theorem for the ordinary cohomology of $\cM_{g,m}$ is Harer’s bound which asserts that for any $m>0$ and $g\ge 1$,
\[
H^{k}(\cM_{g,m}) = 0 \ \ \  \text{for all } k > 4g - 4 + m.
\]

Let us investigate which bound on the non-vanishing cohomology classes follows from the above Corollary (\ref{2: Main equality}). That Corollary says that every non-zero, say of degree $k$, cohomology class in  $H^\bu(\cM_{g,m})$ can be represented by a linear combination of genus $g$ ribbon quivers $\Ga\in \sRGC^{g,m}_{d+1}$ of degree $|\Ga|$ satisfying the equation 
\Beq\label{2: degree p versus degree Ga}
k=6g-6-2dg+3m-dm+d-|\Ga|.
\Eeq
As $\Ga$ is at most 4-valent, we can assume without loss of generality that $\Ga$ has $p_2$ bivalent vertices,  $p_3$ trivalent vertices and $p_4$ four-valent ones. Then 
$$
\# V(\Ga)=p_2+p_3+p_4,\ \ \ \#E(\Ga)=\frac{1}{2}(2p_2 + 3p_3+4p_4)=p_2+ \frac{3}{2}p_3 +2p_4, \ \ \  \#E(\Ga)-\# V(\Ga)=\frac{1}{2}p_3 +p_4.
$$
Its genus is given by the formula
$$
 2g= 2+ \# E(\Ga) - \# V(\Ga)- \# B(\Ga)= 
 \frac{1}{2}p_3+p_4+2-m.
 $$
 implying
 $
 \frac{1}{2}p_3+p_4=2g-2+m
 $.
On the other hand, the cohomological degree of such a graph $\Ga$ is given by
\Beqrn
|\Ga| &=& (d+1)(\#V(\Ga)-1) + d\#E(\Ga)\\
      &=& (d+1)(p_2+p_3+p_4-1) -d(p_2+ \frac{3}{2}p_3+2p_4) \\
      %&=& p_2  +p_3(d+1 - \frac{3}{2}d ) +p_4(d+1-2d) -d-1 \\
      %&=& p_2  + \frac{1}{2}p_3(2-d) + p_4(1-d) -d-1 \\
      &=& p_2 + \frac{1}{2}p_3 + (1-d)(2g-2+ m)  -d-1
      %&=& p_2 + \frac{1}{2}p_3  -2  \ \ \text{for $d=1$}\\ 
\Eeqrn
Hence the above formula (\ref{2: degree p versus degree Ga}) for the degree $k$ of the associated to $\Ga$ cohomology class takes the form
\Beqrn
k%&=&-(p_2 +\frac{1}{2}p_3)+  6g-6-2dg+3m-dm+d-2g+2- m +2gd-2d+dm +d+1\\
&=& -(p_2 +\frac{1}{2}p_3) + 4g+2m-3  \ \ \text{or, equivalently,}\\
&=& %p_4-p_2 -2g-m+2   + 4g-3+2m
p_4-p_2+2g+m-1.
\Eeqrn
Note that the final answer does not depend on the integer parameter $d$ (as expected). As $p_2\geq 2$ and $p_3\geq 0$, we conclude that the cohomology groups  $H^k(\cM_{g,m})$ can be non-zero only in the range
$$
k\leq  + 4g+2m-5=4g-4+m + (m-1)
$$
which does {\it not}\, contradict  Harer's vanishing bound. Cycles in 
$\sRGC^{g,m}_{d+1}$ which are at most trivalent can generate cohomology classes only in degrees
$$
k=-p_2+2g+m-1\leq 2g+m-3.
$$
This calculation  confirms again the  observation that the {\it at most four}\, valency condition on the generators of a ``small" ribbon graph complex 
is sharp (it is well-known that $H^{4g-4+m}(\cM_{g,m})\neq 0$ for $g\geq 1, m\geq 1$).

%{\bf Check: Assume $d=0$}. Then 
%$$
%|\Ga|=p_2+ p_3 +p_4=p_2+ \frac{1}{2}p_3 + 2g-3+m
%$$
%and
%$$
%k=6g-6+3m -|\Ga|=-(p_2+ \frac{1}{2}p_3) + 4g-3+2m=p_4-p_2+ %4g-3+2m-2g+2-m=p_4-p_2+2g+m-1.
%$$
\subsection{Example} By Corollary {\ref{2: Corollary on ORCG to dRGC}}, 
every cycle representative $\Ga$ of a cohomology class from $H^\bu(\ORGC_{d+1})$  must be a coboundary in the complex $\dRGC_{d+1}$, i.e.
$$
\Ga=\delta \theta
$$
for some linear combination of ribbon graphs $\theta\in \dRGC_{d+1}$ which {\it have at least one closed path of directed edges}. This gives a strategy on how to find such classes explicitly. For example,
it is a straightforward calculation to check the equality
$$
\delta \left( 
 \Ba{c}\resizebox{13mm}{!}{
\xy
   {\ar@{<-}@/^0.4pc/(5,0)*{};(-5,0)*{\bu}};
 {\ar@{<-}@/^0.4pc/(-5,0)*{\bu};(5,0)*{\bu}};
 (0,0)*{_{\bar{2}}},
 (0,-4)*{_{\bar{3}}},
(0,4)*{_{\bar{1}}},
 (-5,0)*{}="a",
(5,0)*{}="b",
(0,-9)*{\bu}="c",
\ar @{<-} "c";"a"
\ar @{<-} "c";"b"
\endxy}\Ea   
-
\Ba{c}\resizebox{13mm}{!}{
\xy
 (0,0)*{_{\bar{2}}},
 (0,-4)*{_{\bar{3}}},
(0,4)*{_{\bar{1}}},
   {\ar@{<-}@/^0.4pc/(5,0)*{};(-5,0)*{\bu}};
 {\ar@{<-}@/^0.4pc/(-5,0)*{\bu};(5,0)*{\bu}};
 (-5,0)*{}="a",
(5,0)*{}="b",
(0,-9)*{\bu}="c",
\ar @{->} "c";"a"
\ar @{->} "c";"b"
\endxy}\Ea    
 \right)
 %%%%%%%%%%%%%%%%%%%%%%
 =
 \Ba{c}\resizebox{11mm}{!}{\xy
 (0,-4)*{_{\bar{2}}},
 (5,-1)*{_{\bar{3}}},
(0,-11)*{_{\bar{1}}},
(0,2)*{\bullet}="1",
(-6,-6)*{\bullet}="2",
(6,-10)*{\bullet}="3",
(0,-18)*{\bullet}="4",
\ar @{<-} "4";"3" <0pt>
\ar @{<-} "4";"2" <0pt>
\ar @{<-} "3";"2" <0pt>
\ar @{<-} "2";"1" <0pt>
\ar @{<-} "3";"1" <0pt>
\endxy}\Ea
-
 \Ba{c}\resizebox{11mm}{!}{\xy
 (0,-4)*{_{\bar{1}}},
 (5,-1)*{_{\bar{3}}},
(0,-11)*{_{\bar{2}}},
(0,2)*{\bullet}="1",
(-6,-10)*{\bullet}="2",
(6,-6)*{\bullet}="3",
(0,-18)*{\bullet}="4",
\ar @{<-} "4";"3" <0pt>
\ar @{<-} "4";"2" <0pt>
\ar @{->} "3";"2" <0pt>
\ar @{<-} "2";"1" <0pt>
\ar @{<-} "3";"1" <0pt>
\endxy}\Ea
-
 \Ba{c}\resizebox{11mm}{!}{\xy
 (0,-4)*{_{\bar{2}}},
 (5,-1)*{_{\bar{1}}},
(0,-11)*{_{\bar{3}}},
(0,2)*{\bullet}="1",
(-6,-6)*{\bullet}="2",
(6,-10)*{\bullet}="3",
(0,-18)*{\bullet}="4",
\ar @{<-} "4";"3" <0pt>
\ar @{<-} "4";"2" <0pt>
\ar @{<-} "3";"2" <0pt>
\ar @{<-} "2";"1" <0pt>
\ar @{<-} "3";"1" <0pt>
\endxy}\Ea
+
 \Ba{c}\resizebox{11mm}{!}{\xy
 (0,-4)*{_{\bar{3}}},
 (5,-1)*{_{\bar{1}}},
(0,-11)*{_{\bar{2}}},
(0,2)*{\bullet}="1",
(-6,-10)*{\bullet}="2",
(6,-6)*{\bullet}="3",
(0,-18)*{\bullet}="4",
\ar @{<-} "4";"3" <0pt>
\ar @{<-} "4";"2" <0pt>
\ar @{->} "3";"2" <0pt>
\ar @{<-} "2";"1" <0pt>
\ar @{<-} "3";"1" <0pt>
\endxy}\Ea=:\ga_{0,3}
$$
The vertices of the linear combination $\Ga$ of ribbon quivers graphs on the right hand side are ordered from the top to the bottom (assuming $d$ is even) while their boundaries are distinguished by symbols $\bar{1},\bar{2},\bar{3}$. This linear combination $\ga_{0,3}$ is a cycle in $\sRGC_{1}$ representing a generator of $H^0(\cM_{0,3})\simeq \K$.  Clearly, the relabelling of the boundaries in $\ga_{0,3}$  using the premuation $(13)\in \bS_3$ gives the equality
$$
(13)\ga_{0,3}=-\ga_{0,3}.
$$
It is a straightforward calculation to check that the ``2+3" relations (\ref{2: 2+3 relations}) in 
$\sRGC_{1}$ imply the following equality of ribbon quivers,
$$
 \Ba{c}\resizebox{11mm}{!}{\xy
 (0,-4)*{_{\bar{2}}},
 (5,-1)*{_{\bar{3}}},
(0,-11)*{_{\bar{1}}},
(0,2)*{\bullet}="1",
(-6,-6)*{\bullet}="2",
(6,-10)*{\bullet}="3",
(0,-18)*{\bullet}="4",
\ar @{<-} "4";"3" <0pt>
\ar @{<-} "4";"2" <0pt>
\ar @{<-} "3";"2" <0pt>
\ar @{<-} "2";"1" <0pt>
\ar @{<-} "3";"1" <0pt>
\endxy}\Ea
- \Ba{c}\resizebox{11mm}{!}{\xy
 (0,-4)*{_{\bar{1}}},
 (5,-1)*{_{\bar{3}}},
(0,-11)*{_{\bar{2}}},
(0,2)*{\bullet}="1",
(-6,-10)*{\bullet}="2",
(6,-6)*{\bullet}="3",
(0,-18)*{\bullet}="4",
\ar @{<-} "4";"3" <0pt>
\ar @{<-} "4";"2" <0pt>
\ar @{->} "3";"2" <0pt>
\ar @{<-} "2";"1" <0pt>
\ar @{<-} "3";"1" <0pt>
\endxy}\Ea
\simeq 
\Ba{c}\resizebox{11mm}{!}{\xy
 (0,-4)*{_{\bar{3}}},
 (5,-1)*{_{\bar{1}}},
(0,-11)*{_{\bar{2}}},
(0,2)*{\bullet}="1",
(-6,-6)*{\bullet}="2",
(6,-10)*{\bullet}="3",
(0,-18)*{\bullet}="4",
\ar @{<-} "4";"3" <0pt>
\ar @{<-} "4";"2" <0pt>
\ar @{<-} "3";"2" <0pt>
\ar @{<-} "2";"1" <0pt>
\ar @{<-} "3";"1" <0pt>
\endxy}\Ea
-
\Ba{c}\resizebox{11mm}{!}{\xy
 (0,-4)*{_{\bar{2}}},
 (5,-1)*{_{\bar{1}}},
(0,-11)*{_{\bar{3}}},
(0,2)*{\bullet}="1",
(-6,-10)*{\bullet}="2",
(6,-6)*{\bullet}="3",
(0,-18)*{\bullet}="4",
\ar @{<-} "4";"3" <0pt>
\ar @{<-} "4";"2" <0pt>
\ar @{->} "3";"2" <0pt>
\ar @{<-} "2";"1" <0pt>
\ar @{<-} "3";"1" <0pt>
\endxy}\Ea
$$
which in turn implies the following equivalent  
representation for the above cycle $\ga_{0,3}$,
$$
\ga_{0,3}= \Ba{c}\resizebox{11mm}{!}{\xy
 (0,-4)*{_{\bar{3}}},
 (5,-1)*{_{\bar{1}}},
(0,-11)*{_{\bar{2}}},
(0,2)*{\bullet}="1",
(-6,-6)*{\bullet}="2",
(6,-10)*{\bullet}="3",
(0,-18)*{\bullet}="4",
\ar @{<-} "4";"3" <0pt>
\ar @{<-} "4";"2" <0pt>
\ar @{<-} "3";"2" <0pt>
\ar @{<-} "2";"1" <0pt>
\ar @{<-} "3";"1" <0pt>
\endxy}\Ea
-
\Ba{c}\resizebox{11mm}{!}{\xy
 (0,-4)*{_{\bar{2}}},
 (5,-1)*{_{\bar{1}}},
(0,-11)*{_{\bar{3}}},
(0,2)*{\bullet}="1",
(-6,-10)*{\bullet}="2",
(6,-6)*{\bullet}="3",
(0,-18)*{\bullet}="4",
\ar @{<-} "4";"3" <0pt>
\ar @{<-} "4";"2" <0pt>
\ar @{->} "3";"2" <0pt>
\ar @{<-} "2";"1" <0pt>
\ar @{<-} "3";"1" <0pt>
\endxy}\Ea
-
 \Ba{c}\resizebox{11mm}{!}{\xy
 (0,-4)*{_{\bar{2}}},
 (5,-1)*{_{\bar{1}}},
(0,-11)*{_{\bar{3}}},
(0,2)*{\bullet}="1",
(-6,-6)*{\bullet}="2",
(6,-10)*{\bullet}="3",
(0,-18)*{\bullet}="4",
\ar @{<-} "4";"3" <0pt>
\ar @{<-} "4";"2" <0pt>
\ar @{<-} "3";"2" <0pt>
\ar @{<-} "2";"1" <0pt>
\ar @{<-} "3";"1" <0pt>
\endxy}\Ea
+
 \Ba{c}\resizebox{11mm}{!}{\xy
 (0,-4)*{_{\bar{3}}},
 (5,-1)*{_{\bar{1}}},
(0,-11)*{_{\bar{2}}},
(0,2)*{\bullet}="1",
(-6,-10)*{\bullet}="2",
(6,-6)*{\bullet}="3",
(0,-18)*{\bullet}="4",
\ar @{<-} "4";"3" <0pt>
\ar @{<-} "4";"2" <0pt>
\ar @{->} "3";"2" <0pt>
\ar @{<-} "2";"1" <0pt>
\ar @{<-} "3";"1" <0pt>
\endxy}\Ea.
$$
This representation makes it evident that
$$
(23)\ga_{0,3}=-\ga_{0,3}.
$$
so that we obtain an isomorphism of $\bS_3$-modules
$$
H^3(\sRGC_1^{0,3})= \text{span}_\K \langle [\ga_{0,3}] \rangle = H^0(\cM_{0,3})\ot \sgn_3=\sgn_3,
$$
in accordance with what the general theory suggest (see 
\S {\ref{ 2: subsec on S_m}}).

\sip

The involution (\ref{2: involution iota}) acts on the above graphs (for $d$ even) as follows
$$
\iota\left(  \Ba{c}\resizebox{11mm}{!}{\xy
 (0,-4)*{_{\bar{2}}},
 (5,-1)*{_{\bar{3}}},
(0,-11)*{_{\bar{1}}},
(0,2)*{\bullet}="1",
(-6,-6)*{\bullet}="2",
(6,-10)*{\bullet}="3",
(0,-18)*{\bullet}="4",
\ar @{<-} "4";"3" <0pt>
\ar @{<-} "4";"2" <0pt>
\ar @{<-} "3";"2" <0pt>
\ar @{<-} "2";"1" <0pt>
\ar @{<-} "3";"1" <0pt>
\endxy}\Ea   \right)
=
- \Ba{c}\resizebox{11mm}{!}{\xy
 (0,-4)*{_{\bar{1}}},
 (5,-1)*{_{\bar{3}}},
(0,-11)*{_{\bar{2}}},
(0,2)*{\bullet}="1",
(-6,-6)*{\bullet}="2",
(6,-10)*{\bullet}="3",
(0,-18)*{\bullet}="4",
\ar @{<-} "4";"3" <0pt>
\ar @{<-} "4";"2" <0pt>
\ar @{<-} "3";"2" <0pt>
\ar @{<-} "2";"1" <0pt>
\ar @{<-} "3";"1" <0pt>
\endxy}\Ea, 
\ \ \ 
\iota\left(
 \Ba{c}\resizebox{11mm}{!}{\xy
 (0,-4)*{_{\bar{1}}},
 (5,-1)*{_{\bar{3}}},
(0,-11)*{_{\bar{2}}},
(0,2)*{\bullet}="1",
(-6,-10)*{\bullet}="2",
(6,-6)*{\bullet}="3",
(0,-18)*{\bullet}="4",
\ar @{<-} "4";"3" <0pt>
\ar @{<-} "4";"2" <0pt>
\ar @{->} "3";"2" <0pt>
\ar @{<-} "2";"1" <0pt>
\ar @{<-} "3";"1" <0pt>
\endxy}\Ea\right)
=
- \Ba{c}\resizebox{11mm}{!}{\xy
 (0,-4)*{_{\bar{2}}},
 (5,-1)*{_{\bar{3}}},
(0,-11)*{_{\bar{1}}},
(0,2)*{\bullet}="1",
(-6,-10)*{\bullet}="2",
(6,-6)*{\bullet}="3",
(0,-18)*{\bullet}="4",
\ar @{<-} "4";"3" <0pt>
\ar @{<-} "4";"2" <0pt>
\ar @{->} "3";"2" <0pt>
\ar @{<-} "2";"1" <0pt>
\ar @{<-} "3";"1" <0pt>
\endxy}\Ea.
$$
Hence
$$
\iota(\ga_{0,3})= \ga_{0.3},
$$
i.e.\  the cohomologically non-trivial cycle $\ga_{0,3}$ belongs to $\sRGC_1^+$ (as expected from 
Corollary {\ref{2: iso of sRGC+,g,m with M}}).

\bip

\bip

{\Large
\section{\bf A proof of Main Theorem {\ref{2: Main theorem}}}
}

\bip

\subsection{Plan} As complexes $\ORGC_{d+1}$ for various $d$ are all isomorphic to each other
(up to degree shifts), it is enough to prove the theorem for $d=0$, when the sign rules are simplest possible; in this case vertices of ribbon quivers have degree $+1$, edges have degree zero, and their orientations reduce to the orderings of vertices (which we assume in the pictures below to run from the top to the bottom).

\sip

The main idea of the proof is to use appropriate filtrations of the complex $\ORGC_1$ such that on the first two pages of the associated spectral sequences the problem reduces essentially to the study of the tensor products of four complexes: 
\Bi
\item[(i)] 
 the dg properad $\mathcal{IB}_\infty$ of strongly homotopy infinitesimal bialgebras whose cohomology was computed in \cite{A}, 
 \item[(ii)] the dg properad $\caD\Lie_\infty$ of strongly homotopy double Lie algebras whose cohomology was computed in \cite{L},
 \item[(iii)] an auxiliary dg free right module $\cA ss_\infty^{\wedge,+}$ over the dg operad  $\cA ss_\infty$ of $A_\infty$-algebras whose cohomology is computed in \S {\ref{3: subsec on plus extensions}} below;
 \item[(iv)] an auxiliary dg free right module $\cA ss_\infty^{\wedge,cyc}$ over the operad  $\cA ss_\infty$  whose cohomology is computed in \S {\ref{3: subsec on Ass_wedge_cyc}} below.   
 \Ei
So we start our proof with reminders about the known cohomology groups of complexes $\cA ss_\infty$,  $\mathcal{IB}_\infty$ and $\caD\Lie_\infty$ (essentially to fix the notation); then we define and compute cohomologies of auxiliary complexes $\cA ss_\infty^{\wedge,+}$ and $\cA ss_\infty^{\wedge,cyc}$. After that we apply this set of  results to show that the cohomology groups $H^\bu(\ORGC_1)$ and $H^\bu(\sRGC_1)$ are isomorphic.

\mip

\subsection{The complex $\cA ss_\infty$} The dg free operad 
$\cA ss_\infty=\{\cA ss_\infty(n)\}_{n\geq 2}$ is generated by degree $+1$ planar corollas with $n\geq 2$ incoming legs (which are totally ordered from the left to the right)
$$
  \underbrace{\Ba{c}\resizebox{20mm}{!}
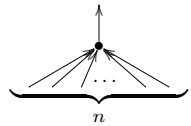
{\xy
(1,-6)*{\ldots},
(-13,-7)*{},
(-8,-7)*{},
(-3,-7)*{},
(7,-7)*{},
(13,-7)*{},
 (0,0)*{\bu}="a",
(0,7)*{}="0",
(-12,-7)*{}="b_1",
(-8,-7)*{}="b_2",
(-3,-7)*{}="b_3",
(8,-7)*{}="b_4",
(12,-7)*{}="b_5",
\ar @{->} "a";"0" <0pt>
\ar @{<-} "a";"b_2" <0pt>
\ar @{<-} "a";"b_3" <0pt>
\ar @{<-} "a";"b_1" <0pt>
\ar @{<-} "a";"b_4" <0pt>
\ar @{<-} "a";"b_5" <0pt>
\endxy}\Ea}_n
$$

and one outgoing leg. Notice that the number of directed paths from the in-edges to the out-edges is equal to $n$. The differential is given on the generators by substituting into the vertex the graph $\Ba{c}\xy
 (0,0)*{\bu}="a",
(0,5)*{\bu}="b",
\ar @{->} "a";"b" <0pt>
\endxy\Ea$ and reattaching the edges among the two newly created vertices in such a way that the total order of edges from the left to right is preserved. Thus
\Beq\label{3: d in Ass_infty}
\delta
  \underbrace{\Ba{c}\resizebox{20mm}{!}{\xy
(1,-6)*{\ldots},
(-13,-7)*{},
(-8,-7)*{},
(-3,-7)*{},
(7,-7)*{},
(13,-7)*{},
 (0,0)*{\bu}="a",
(0,7)*{}="0",
(-12,-7)*{}="b_1",
(-8,-7)*{}="b_2",
(-3,-7)*{}="b_3",
(8,-7)*{}="b_4",
(12,-7)*{}="b_5",
\ar @{->} "a";"0" <0pt>
\ar @{<-} "a";"b_2" <0pt>
\ar @{<-} "a";"b_3" <0pt>
\ar @{<-} "a";"b_1" <0pt>
\ar @{<-} "a";"b_4" <0pt>
\ar @{<-} "a";"b_5" <0pt>
\endxy}\Ea}_n
=
\sum_{A\subsetneq [n]} 
\Ba{c}\resizebox{22mm}{!}{\xy
(1.6,-7)*{...},
(-13,-7)*{},
(-8,-7)*{},
(-3,-7)*{},
(7,-7)*{},
(13,-7)*{},
 (0,0)*{\bu}="a",
(0,7)*{}="0",
(-12,-7)*{}="b_1",
(-9,-7)*{}="b_2",
(-6,-7)*{...},
(-3,-7)*{\bu}="b_3",
(-9,-14)*{}="c1",
(-5,-14)*{}="c2",
(-1,-14)*{...},
(3,-14)*{}="c3",
(7,-7)*{}="b_4",
(11,-7)*{}="b_5",
(-3,-17)*{\underbrace{\ \ \ \ \ \ \ \ \ \ \ \ }_A},
\ar @{->} "a";"0" <0pt>
\ar @{<-} "a";"b_2" <0pt>
\ar @{<-} "a";"b_3" <0pt>
\ar @{<-} "a";"b_1" <0pt>
\ar @{<-} "a";"b_4" <0pt>
\ar @{<-} "a";"b_5" <0pt>
\ar @{->} "c1";"b_3" <0pt>
\ar @{->} "c2";"b_3" <0pt>
\ar @{->} "c3";"b_3" <0pt>
\endxy}\Ea
\Eeq
where the summation runs over connected proper subsets $A$ (of cardinality $\geq 2$) of the totally ordered set $[n]$.
Its cohomology $H^\bu(\cA ss_\infty)$ was proven in \cite{GK} to be the operad $\cA ss$ generated by the planar corolla $
  \Ba{c}\resizebox{8mm}{!}{  \xy
(0,6)*{}="1";
    (0,0.2)*{\bu}="L";
  (-4,-5)*{}="C";
   (+4,-5)*{}="D";
\ar @{->} "D";"L" <0pt>
\ar @{->} "C";"L" <0pt>
\ar @{<-} "1";"L" <0pt>
 \endxy}
 \Ea
$ (whose two incoming legs are totally ordered from the left to the right) modulo the associativity relation
 \Beq\label{3: Ass operad relation}
  \Ba{c}\resizebox{13.5mm}{!}{  \xy
(-9,10)*{^0}="1";
    (-9,+3)*{\bbu}="L";
 (-14,-3.5)*{\bbu}="B";
 (-20,-12)*+{_1}="b1";
 (-8,-12)*+{_2}="b2";
  (-3,-4)*{_3}="C";
\ar @{->} "C";"L" <0pt>
\ar @{->} "B";"L" <0pt>
\ar @{<-} "B";"b1" <0pt>
\ar @{<-} "B";"b2" <0pt>
\ar @{<-} "1";"L" <0pt>
 \endxy}
 \Ea
 +
  \Ba{c}\resizebox{13.5mm}{!}{  \xy
(9,10)*{^0}="1";
    (9,+3)*{\bbu}="L";
 (14,-3.5)*{\bbu}="B";
 (20,-12)*+{_3}="b1";
 (8,-12)*+{2}="b2";
  (3,-4)*{1}="C";
\ar @{->} "C";"L" <0pt>
\ar @{->} "B";"L" <0pt>
\ar @{<-} "B";"b1" <0pt>
\ar @{<-} "B";"b2" <0pt>
\ar @{<-} "1";"L" <0pt>
 \endxy}
 \Ea=0.
 \Eeq
 Notice that these relations are given precisely  by those ribbon graphs in (\ref{2: twisted ass-type relations}) which have the maximal possible number of paths from the input legs to the unique outgoing edge. We shall use this fact below.
 
 \sip
 
 Denote by $\cA ss_\infty^{op}$ and $\cA ss^{op}$ the versions of the above two operads in which the directions of all arrows are reversed.

 \sip

 \subsection{An auxiliary ``plus" module}\label{3: subsec on plus extensions}
   Let $\cA ss_\infty^{\wedge, +}=\{\cA ss_\infty^{\wedge, +}(n)\}_{n\geq 1}$  be a dg free right module over the dg operad $\cA ss_\infty$ generated by degree +1 planar corollas of the form
 $$
  \underbrace{\Ba{c}\resizebox{20mm}{!}
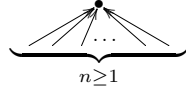
{\xy
(1,-6)*{\ldots},
(-13,-7)*{},
(-8,-7)*{},
(-3,-7)*{},
(7,-7)*{},
(13,-7)*{},
 (0,0)*{\bu}="a",
(0,7)*{}="0",
(-12,-7)*{}="b_1",
(-8,-7)*{}="b_2",
(-3,-7)*{}="b_3",
(8,-7)*{}="b_4",
(12,-7)*{}="b_5",
%
%\ar @{->} "a";"0" <0pt>
\ar @{<-} "a";"b_2" <0pt>
\ar @{<-} "a";"b_3" <0pt>
\ar @{<-} "a";"b_1" <0pt>
\ar @{<-} "a";"b_4" <0pt>
\ar @{<-} "a";"b_5" <0pt>
\endxy}\Ea}_{n\geq 1}
$$
The adjective {\it planar}\, means that the incoming edges are totally ordered from the left to the right.
The differential $\delta^+$ is given on elements $\Ga$ by the standard sum over vertices 
$$
\delta^+\Ga=\sum_{v\in V(\Ga)}\delta^+_v \Ga
$$
 where $\delta_v^+$ is given on $\cA ss_\infty$-vertices by (\ref{3: d in Ass_infty}) and on vertices $v$ with no outgoing edge (i.e.\ on targets) as follows, 
\Beq\label{3: d in Ass_wedge_plus_infty}
\delta^+ \hspace{-2mm}  \Ba{c}\resizebox{2.5mm}{!}{  \xy
    (0,2)*{\bu}="L";
  (0,-4)*{}="C";
\ar @{->} "C";"L" <0pt>
 \endxy}
 \Ea=0, \ \ \
 \delta^+ \hspace{-3mm}   \Ba{c}\resizebox{8mm}{!}{  \xy
(0,6)*{}="1";
    (0,0.2)*{\bu}="L";
  (-4,-5)*{}="C";
   (+4,-5)*{}="D";
\ar @{->} "D";"L" <0pt>
\ar @{->} "C";"L" <0pt>
%\ar @{<-} "1";"L" <0pt>
%
 \endxy}
 \Ea
 =  \Ba{c}\resizebox{8mm}{!}{  \xy
(0,6)*{\bu}="1";
    (0,0.2)*{\bu}="L";
  (-4,-5)*{}="C";
   (+4,-5)*{}="D";
\ar @{->} "D";"L" <0pt>
\ar @{->} "C";"L" <0pt>
\ar @{<-} "1";"L" <0pt>
 \endxy}
 \Ea, \ \ \ \ \
 \delta^+\hspace{-3mm}
  \underbrace{\Ba{c}\resizebox{20mm}{!}{\xy
(1,-6)*{\ldots},
(-13,-7)*{},
(-8,-7)*{},
(-3,-7)*{},
(7,-7)*{},
(13,-7)*{},
 (0,0)*{\bu}="a",
(0,7)*{}="0",
(-12,-7)*{}="b_1",
(-8,-7)*{}="b_2",
(-3,-7)*{}="b_3",
(8,-7)*{}="b_4",
(12,-7)*{}="b_5",
%
%\ar @{->} "a";"0" <0pt>
\ar @{<-} "a";"b_2" <0pt>
\ar @{<-} "a";"b_3" <0pt>
\ar @{<-} "a";"b_1" <0pt>
\ar @{<-} "a";"b_4" <0pt>
\ar @{<-} "a";"b_5" <0pt>
\endxy}\Ea}_{n\geq 3}
=
\Ba{c}\resizebox{20mm}{!}{\xy
(1,-6)*{\ldots},
(-13,-7)*{},
(-8,-7)*{},
(-3,-7)*{},
(7,-7)*{},
(13,-7)*{},
 (0,0)*{\bu}="a",
(0,7)*{\bu}="0",
(-12,-7)*{}="b_1",
(-8,-7)*{}="b_2",
(-3,-7)*{}="b_3",
(8,-7)*{}="b_4",
(12,-7)*{}="b_5",
\ar @{->} "a";"0" <0pt>
\ar @{<-} "a";"b_2" <0pt>
\ar @{<-} "a";"b_3" <0pt>
\ar @{<-} "a";"b_1" <0pt>
\ar @{<-} "a";"b_4" <0pt>
\ar @{<-} "a";"b_5" <0pt>
\endxy}\Ea\hspace{-3mm}
%%%%%%%%%%
+
\sum_{A\subsetneq [n]\atop \# A\geq 2} 
\Ba{c}\resizebox{22mm}{!}{\xy
(1.6,-7)*{...},
(-13,-7)*{},
(-8,-7)*{},
(-3,-7)*{},
(7,-7)*{},
(13,-7)*{},
 (0,0)*{\bu}="a",
(0,7)*{}="0",
(-12,-7)*{}="b_1",
(-9,-7)*{}="b_2",
(-6,-7)*{...},
(-3,-7)*{\bu}="b_3",
(-9,-14)*{}="c1",
(-5,-14)*{}="c2",
(-1,-14)*{...},
(3,-14)*{}="c3",
(7,-7)*{}="b_4",
(11,-7)*{}="b_5",
(-3,-17)*{\underbrace{\ \ \ \ \ \ \ \ \ \ \ \ }_A},
%
%\ar @{->} "a";"0" <0pt>
\ar @{<-} "a";"b_2" <0pt>
\ar @{<-} "a";"b_3" <0pt>
\ar @{<-} "a";"b_1" <0pt>
\ar @{<-} "a";"b_4" <0pt>
\ar @{<-} "a";"b_5" <0pt>
\ar @{->} "c1";"b_3" <0pt>
\ar @{->} "c2";"b_3" <0pt>
\ar @{->} "c3";"b_3" <0pt>
\endxy}\Ea.
\Eeq
Using a filtration of each complex $\cA ss_\infty^{\wedge,+}(n)$, $n\geq 1$, by the number of univalent vertices in it is easy to conclude that
\Beq\label{3: H(Ass_infty_wedge+plus)}
H^\bu(\cA ss_\infty^{\wedge,+}(1))=\K[-1], \ \ \ \ H^\bu(\cA ss_\infty^{\wedge,+}(n))=0\ \ 
\text{for 
$n\geq 2$},
\Eeq
i.e.\ the cohomology of $\cA ss_\infty^{\wedge,+}$ is 1-dimensional and is spanned by \hspace{-2mm}
$ \Ba{c}\resizebox{2.5mm}{!}{  \xy
    (0,2)*{\bu}="L";
  (0,-4)*{}="C";
\ar @{->} "C";"L" <0pt>
 \endxy}
 \Ea$.

\sip

We denote by $\cA ss_\infty^{\vee,+}$ a copy of $\cA ss_\infty^{\wedge,+}$ 
in which directions of all arrows in the generators are reversed.

\subsection{An auxiliary complex $\Ass_\infty^{\wedge,cyc}$}\label{3: subsec on Ass_wedge_cyc}  Let $\cA ss_\infty^{\wedge, cyc}=\{\cA ss_\infty^{\wedge,cyc}(n)\}_{n\geq 2}$ be a dg free right module over the dg operad $\cA ss_\infty$ generated by degree +1 ribbon corollas of the form
 $$
\Ba{c}\resizebox{20mm}{!}{\xy
(2,-6)*{\ldots},
(-13,+7)*{},
(-8,-1)*{},
(-3,-7)*{},
(7,-7)*{},
(13,-7)*{},
 (0,0)*{\circ}="a",
(0,7)*{}="0",
(-7,+7)*+{_k}="b_1",
(-9,-3)*+{_1}="b_2",
(-3,-9)*+{_2}="b_3",
(8,-7)*+{}="b_4",
(5,9)*+{_{k-1}}="b_5",
%
%\ar @{->} "a";"0" <0pt>
\ar @{<-} "a";"b_2" <0pt>
\ar @{<-} "a";"b_3" <0pt>
\ar @{<-} "a";"b_1" <0pt>
\ar @{<-} "a";"b_4" <0pt>
\ar @{<-} "a";"b_5" <0pt>
\endxy}\Ea, \ \ \ {k\geq 2}.
$$
 The adjective {\it ribbon}\, means that the incoming edges are {\it cyclically ordered}, say anticlockwise. The cyclic subgroup $\Z_k\subset \bS_k$ generated by the permutation $(1,2,...,k)$ acts trivially on such a generator.
The finite set   $\{12...n\}$ equipped with the cyclic ordering  $1<2<...<n<1$ is denoted by $((n))$.

\sip
 
 The differential in $\cA ss_\infty^{\wedge, cyc}$ is defined by
\Beq\label{3: d in Ass_wedge_infty}
\delta
\Ba{c}\resizebox{20mm}{!}{\xy
(2,-6)*{\ldots},
(-13,+7)*{},
(-8,-1)*{},
(-3,-7)*{},
(7,-7)*{},
(13,-7)*{},
 (0,0)*{\circ}="a",
(0,7)*{}="0",
(-7,+7)*+{_n}="b_1",
(-9,-3)*+{_1}="b_2",
(-3,-9)*+{_2}="b_3",
(8,-7)*+{}="b_4",
(6,7)*+{}="b_5",
%
%\ar @{->} "a";"0" <0pt>
\ar @{<-} "a";"b_2" <0pt>
\ar @{<-} "a";"b_3" <0pt>
\ar @{<-} "a";"b_1" <0pt>
\ar @{<-} "a";"b_4" <0pt>
\ar @{<-} "a";"b_5" <0pt>
\endxy}\Ea
=
\sum_{A\subsetneq [n]\atop \# A\geq 2} 
\Ba{c}\resizebox{20mm}{!}{\xy
%(1.6,-7)*{...},
(-13,-7)*{},
(-8,-7)*{},
(-3,-7)*{},
(7,-7)*{},
(13,-7)*{},
 (0,0)*{\circ}="a",
(0,7)*{}="0",
(-6,8)*{}="b_1",
(-9,-5)*{}="b_2",
%(-6,-7)*{...},
(-3,-7)*{\circ}="b_3",
(-9,-14)*{}="c1",
(-5,-14)*{}="c2",
(-1,-14)*{...},
(3,-14)*{}="c3",
(7,-5)*{}="b_4",
(6,8)*{}="b_5",
(-3,-17)*{\underbrace{\ \ \ \ \ \ \ \ \ \ \ \ }_A},
%
%\ar @{->} "a";"0" <0pt>
\ar @{<-} "a";"b_2" <0pt>
\ar @{<-} "a";"b_3" <0pt>
\ar @{<-} "a";"b_1" <0pt>
\ar @{<-} "a";"b_4" <0pt>
\ar @{<-} "a";"b_5" <0pt>
\ar @{->} "c1";"b_3" <0pt>
\ar @{->} "c2";"b_3" <0pt>
\ar @{->} "c3";"b_3" <0pt>
\endxy}\Ea,
\Eeq
where the summation runs over proper connected subsets $A$ of $((n))$ of cardinality $\geq 2$. Note that each such subset $A$ comes equipped with the induced total order of incoming edges. The right hand side is understood as a sum of ribbon graphs with labelled hairs. For example
$$
\delta \Ba{c}\resizebox{16mm}{!}{ \xy
 (-7,2)*+{_1}="1",
 (0,0)*+{_2}="2",
 (7,2)*+{_3}="3",
(0,8)*{\circ}="c",
\ar @{<-} "c";"1" <0pt>
\ar @{<-} "c";"2" <0pt>
\ar @{<-} "c";"3" <0pt>
\endxy}\Ea
=
  \Ba{c}\resizebox{13.5mm}{!}{  \xy
    (-9,+3)*{\circ}="L";
 (-14,-3.5)*{\circ}="B";
 (-20,-12)*+{_1}="b1";
 (-8,-12)*+{_2}="b2";
  (-3,-4)*{_3}="C";
\ar @{->} "C";"L" <0pt>
\ar @{->} "B";"L" <0pt>
\ar @{<-} "B";"b1" <0pt>
\ar @{<-} "B";"b2" <0pt>
 \endxy}
 \Ea
 +
  \Ba{c}\resizebox{13.5mm}{!}{  \xy
    (-9,+3)*{\circ}="L";
 (-14,-3.5)*{\circ}="B";
 (-20,-12)*+{_2}="b1";
 (-8,-12)*+{_3}="b2";
  (-3,-4)*{_1}="C";
\ar @{->} "C";"L" <0pt>
\ar @{->} "B";"L" <0pt>
\ar @{<-} "B";"b1" <0pt>
\ar @{<-} "B";"b2" <0pt>
 \endxy}
 \Ea
 +
  \Ba{c}\resizebox{13.5mm}{!}{  \xy
    (-9,+3)*{\circ}="L";
 (-14,-3.5)*{\circ}="B";
 (-20,-12)*+{_3}="b1";
 (-8,-12)*+{_1}="b2";
  (-3,-4)*{_2}="C";
\ar @{->} "C";"L" <0pt>
\ar @{->} "B";"L" <0pt>
\ar @{<-} "B";"b1" <0pt>
\ar @{<-} "B";"b2" <0pt>
 \endxy}
 \Ea
 $$
 
\subsubsection{\bf Lemma} {\it  The cohomology $H^\bu(\cA ss_\infty^{\wedge,cyc})$ equals  a right module 
$\cA ss^{\wedge,cyc}$ over the operad $\cA ss$ which, by definition, is generated by the ribbon corolla
 $
  \Ba{c}\resizebox{8mm}{!}{  \xy
    (0,0.2)*{\circ}="L";
  (-4,-5)*{}="C";
   (+4,-5)*{}="D";
\ar @{->} "D";"L" <0pt>
\ar @{->} "C";"L" <0pt>
 \endxy}
 \Ea
$
 modulo the relation}
 \Beq\label{3: relation for Ass_wedge_cyc}
   \Ba{c}\resizebox{13.5mm}{!}{  \xy
    (-9,+3)*{\circ}="L";
 (-14,-3.5)*{\circ}="B";
 (-20,-12)*+{_1}="b1";
 (-8,-12)*+{_2}="b2";
  (-3,-4)*{_3}="C";
\ar @{->} "C";"L" <0pt>
\ar @{->} "B";"L" <0pt>
\ar @{<-} "B";"b1" <0pt>
\ar @{<-} "B";"b2" <0pt>
 \endxy}
 \Ea
 +
  \Ba{c}\resizebox{13.5mm}{!}{  \xy
    (-9,+3)*{\circ}="L";
 (-14,-3.5)*{\circ}="B";
 (-20,-12)*+{_2}="b1";
 (-8,-12)*+{_3}="b2";
  (-3,-4)*{_1}="C";
\ar @{->} "C";"L" <0pt>
\ar @{->} "B";"L" <0pt>
\ar @{<-} "B";"b1" <0pt>
\ar @{<-} "B";"b2" <0pt>
 \endxy}
 \Ea
 +
  \Ba{c}\resizebox{13.5mm}{!}{  \xy
    (-9,+3)*{\circ}="L";
 (-14,-3.5)*{\circ}="B";
 (-20,-12)*+{_3}="b1";
 (-8,-12)*+{_1}="b2";
  (-3,-4)*{_2}="C";
\ar @{->} "C";"L" <0pt>
\ar @{->} "B";"L" <0pt>
\ar @{<-} "B";"b1" <0pt>
\ar @{<-} "B";"b2" <0pt>
 \endxy}
 \Ea
=0. 
\Eeq
 
 \begin{proof} We have to show that $H^\bu(\cA ss_\infty^{\wedge,cyc}(n))=\cA ss^{\wedge,cyc}(n)$ for all $n\geq 2$. We can assume without loss of generality that the in-legs of the generators of $\cA ss_\infty^{\wedge,cyc}(n)$    are labelled by the set $((n))$ 
 in the standard anticlockwise order. The case $n=2$ is clear, so we assume from now on that $n\geq 3$.
 
 \sip

 As a graded vector space,   $\cA ss_\infty^{\wedge,cyc}(n)$ can be understood as a direct sum of subspaces
 \Beq\label{3: ass_wedge_cyc_decomposition}
 \cA ss_\infty^{\wedge,cyc}(n)=\sum_{k=2}^{n}\sum_{n=p_1+...+p_k} 
 \left(\cA ss^{1+}(p_1)\ot \cA ss^{1+}(p_2)\ot \ldots \ot \cA ss^{1+}(p_k)  \right)_{\Z_k}.
 \Eeq
 where $\cA ss^{1+}(p)$ is a linear subspace of $\cA ss^{\wedge+}(p)$, $p\geq 1$, spanned by elements  with the  target of valency precisely one. 
A generator of the $k$-th summand can be pictorially understood as a cyclically ordered set of $k$ elements of $\cA ss_\infty^{1+}$ whose unique sources are glued into one white vertex making the latter look like  a ribbon corolla  $
\Ba{c}\resizebox{15mm}{!}{\xy
(2,-6)*{\ldots},
(-13,+7)*{},
(-8,-1)*{},
(-3,-7)*{},
(7,-7)*{},
(13,-7)*{},
 (0,0)*{\circ}="a",
(0,7)*{}="0",
(-7,+7)*+{}="b_1",
(-9,-3)*+{}="b_2",
(-3,-9)*+{}="b_3",
(8,-7)*+{}="b_4",
(6,7)*+{}="b_5",
%
%\ar @{->} "a";"0" <0pt>
\ar @{<-} "a";"b_2" <0pt>
\ar @{<-} "a";"b_3" <0pt>
\ar @{<-} "a";"b_1" <0pt>
\ar @{<-} "a";"b_4" <0pt>
\ar @{<-} "a";"b_5" <0pt>
\endxy}\Ea
$ with $k$ edges/legs attached.
The labelled legs of such a generator of $\cA ss_\infty^{\wedge,cyc}(n)$ are hence grouped in totally ordered proper subsets of cardinalities $p_i$, $i\in ((k))$.

\sip
 
There is an obvious epimorphism 
 $$
 s: \cA ss_\infty^{\wedge,cyc}(n) \lon \cA ss^{\wedge,cyc}(n).
 $$
 Consider increasing  filtrations 
 $$
 F_0\subset ...\subset F_p\subset F_{p+1}\subset...
 $$
 of the both sides of this morphism by the number $p$ of vertices in the brunch containing the leg labelled by $1$. Let
  $$
 s_r: \cE_r\cA ss_\infty^{\wedge,cyc}(n)) \lon \cE_r\cA ss^{\wedge,cyc} 
 $$
 be the induced collection of morphisms of $r$-th pages of the associated graded sequences.
  If we show that 
  $$
 s_2: \cE_2\cA ss_\infty^{\wedge,cyc}(n)) \lon \cE_2\cA ss^{\wedge,cyc}(n).
 $$
 is an isomorphism, we are done.
 
 \sip
 
 As $\cA ss^{\wedge,cyc}(n)$ has zero differential, each page  $\cE_r\cA ss^{\wedge,cyc}$ is equal to
  $\cE_0\cA ss^{\wedge,cyc}$, the associated graded of  $\cA ss^{\wedge,cyc}$.
 Using iteratively relations (\ref{3: relation for Ass_wedge_cyc}) and (\ref{3: Ass operad relation}), one concludes easily that $\cE_0\cA ss^{\wedge,cyc}(n))$ for $n\geq 3$ is a vector space of dimension $n-1$ spanned by the following set of basis vectors
 $$
 \left\{   (12...n)^{r}\cdot \Ba{c}\resizebox{17.5mm}{!}{  \xy
    (-9,+3)*{\circ}="L";
 (-14,-3.5)*{\circ}="B";
 (-24,-18)*{\circ}="B2";
 (-20,-12)*+{...}="B1";
 (-15,-19)*{}="b0";
 (-8,-12)*+{_{n-1}}="b2";
 (-29,-25)*+{_1}="b3";
  (-19,-25)*+{_2}="b4";
  (-3,-4)*+{_{n}}="C";
\ar @{->} "C";"L" <0pt>
\ar @{->} "B";"L" <0pt>
\ar @{<-} "B";"B1" <0pt>
\ar @{<-} "B";"b2" <0pt>
\ar @{<-} "B1";"B2" <0pt>
\ar @{<-} "B1";"b0" <0pt>
\ar @{<-} "B2";"b3" <0pt>
\ar @{<-} "B2";"b4" <0pt>
 \endxy}\Ea  \right\}_{r\in \{0,\ldots, n-2\}}
 $$
 where the symbol
 $(12...n)^{r}\cdot $ means the action of the $r$-th power of the cyclic permutation $(12...n)$ on the labels of the legs of the shown graph.
  
  \sip
  
  Assume that in the notation (\ref{3: ass_wedge_cyc_decomposition}) it is the brunch  $\cA ss^{1+}(p_1)$
  which contains the leg labelled by $1$. Then the complex $\cE_0\cA ss^{\wedge,cyc}_\infty(n)$ is isomorphic to the tensor product of the trivial complex $\cA ss^{1+}_\infty(p_1)$ and the complex
   $\cA ss^{\wedge+}_\infty(p_2+...+p_q)$ whose cohomology is non-trivial if and only if $p_2+...+p_q=1$ in which case it equals $\K$.
   Thus the next page $\cE_1\cA ss^{\wedge,cyc}_\infty(n))$   is isomorphic to the direct sum of $n-1$ complexes
   $C_j\simeq \cA ss_\infty(n-1)^{non-sigma}$, $j\in \{2,..., n\}$; here we understand $C_j$ as a dg subspace of $\cA ss_\infty(n-1)$
   spanned by standard planar trees  whose input legs are labelled from the left to the right by the 
   integers $j+1,j+2,..., n,1,2,..., j-1$. Hence the next page $\cE_2\cA ss^{\wedge,cyc}(n)$ is a trivial complex which is identical to $\cE_0\cA ss^{\wedge,cyc}(n))$. Hence the morphism of second pages $s_2$ is an isomorphism. By the comparison theorem of spectral sequences,  the proof of the Lemma is completed.
 \end{proof}

 Denote by $\cA ss_\infty^{\vee,cyc}$ a copy of the dg module $\cA ss_\infty^{\wedge,cyc}$ obtained by reversing directions of arrows of all edges in the generators.

\subsection{Reminder on the dg properad $\mathcal{IB}_\infty$ }\label{3: subsec on IB_infty}
 By definition the properad $\mathcal{IB}_\infty$ of strongly homotopy infinitesimal bialgebras 
is a dg free properad generated by planar $(m,n)$-corollas \cite{A}
\Beq\label{3: (m,n) corolla}
\Ba{c}\resizebox{11mm}{!}{ \xy
(0,8)*{\overbrace{\  \ \ \ \ \ \ \ \ \ }^{m\geq 2}},
(0,-8)*{\underbrace{\  \ \ \ \ \ \ \ \ \ }_{n\geq 2}},
(0,4.5)*+{...},
(0,-4.5)*+{...},
(0,0)*{\bu}="o",
(-5,5)*{}="1",
(-3,5)*{}="2",
(3,5)*{}="3",
(5,5)*{}="4",
(-3,-5)*{}="5",
(3,-5)*{}="6",
(5,-5)*{}="7",
(-5,-5)*{}="8",
\ar @{->} "o";"1" <0pt>
\ar @{->} "o";"2" <0pt>
\ar @{->} "o";"3" <0pt>
\ar @{->} "o";"4" <0pt>
\ar @{<-} "o";"5" <0pt>
\ar @{<-} "o";"6" <0pt>
\ar @{<-} "o";"7" <0pt>
\ar @{<-} "o";"8" <0pt>
\endxy}\Ea
\Eeq
 with $m\geq 1$, $n\geq 1$, $m+n\geq 3$.
Below we understand such a corolla as a vertex of a ribbon graph which 
essentially means that the set of out-legs (respectively in-legs) is totally ordered. 
The differential $\delta$ in $\mathcal{IB}_\infty$ acts on such an $(m,n)$-corolla by substituting into its vertex 
$\bu$ the graph  $\Ba{c}\xy
 (0,0)*{\bu}="a",
(0,5)*{\bu}="b",
\ar @{->} "a";"b" <0pt>
\endxy\Ea$ and reattaching the edges among the two newly created vertices in such a way that the cyclic structure is preserved (which means in this case that total orders are respected).
Thus
\Beq\label{3: d in IB_infty}
\delta
\Ba{c}\resizebox{11mm}{!}{ \xy
(0,8)*{\overbrace{\  \ \ \ \ \ \ \ \ \ }^{m\geq 2}},
(0,-8)*{\underbrace{\  \ \ \ \ \ \ \ \ \ }_{n\geq 2}},
(0,4.5)*+{...},
(0,-4.5)*+{...},
(0,0)*{\bu}="o",
(-5,5)*{}="1",
(-3,5)*{}="2",
(3,5)*{}="3",
(5,5)*{}="4",
(-3,-5)*{}="5",
(3,-5)*{}="6",
(5,-5)*{}="7",
(-5,-5)*{}="8",
\ar @{->} "o";"1" <0pt>
\ar @{->} "o";"2" <0pt>
\ar @{->} "o";"3" <0pt>
\ar @{->} "o";"4" <0pt>
\ar @{<-} "o";"5" <0pt>
\ar @{<-} "o";"6" <0pt>
\ar @{<-} "o";"7" <0pt>
\ar @{<-} "o";"8" <0pt>
\endxy}\Ea
=
\Ba{c}\resizebox{11mm}{!}{ \xy
(0,9.5)*+{...},
(0,-4.5)*+{...},
(0,5)*{\bu}="oo",
(0,0)*{\bu}="o",
(-5,10)*{}="1",
(-3,10)*{}="2",
(3,10)*{}="3",
(5,10)*{}="4",
(-3,-5)*{}="5",
(3,-5)*{}="6",
(5,-5)*{}="7",
(-5,-5)*{}="8",
\ar @{->} "oo";"1" <0pt>
\ar @{->} "oo";"2" <0pt>
\ar @{->} "oo";"3" <0pt>
\ar @{->} "oo";"4" <0pt>
\ar @{->} "o";"oo" <0pt>
\ar @{<-} "o";"5" <0pt>
\ar @{<-} "o";"6" <0pt>
\ar @{<-} "o";"7" <0pt>
\ar @{<-} "o";"8" <0pt>
\endxy}\Ea
+
\sum_{A\subsetneq [n]\atop \# A\geq 2} 
\Ba{c}\resizebox{22mm}{!}{\xy
(1.6,-7)*{...},
(-13,-7)*{},
(-8,-7)*{},
(-3,-7)*{},
(7,-7)*{},
(13,-7)*{},
 (0,0)*{\bu}="a",
 (0,6)*{...},
(-5,7)*{}="01",
(-3,7)*{}="02",
(3,7)*{}="03",
(5,7)*{}="04",
(-12,-7)*{}="b_1",
(-9,-7)*{}="b_2",
(-6,-7)*{...},
(-3,-7)*{\bu}="b_3",
(-9,-14)*{}="c1",
(-5,-14)*{}="c2",
(-1,-14)*{...},
(3,-14)*{}="c3",
(7,-7)*{}="b_4",
(11,-7)*{}="b_5",
(-3,-17)*{\underbrace{\ \ \ \ \ \ \ \ \ \ \ \ }_A},
\ar @{->} "a";"01" <0pt>
\ar @{->} "a";"02" <0pt>
\ar @{->} "a";"03" <0pt>
\ar @{->} "a";"04" <0pt>
\ar @{<-} "a";"b_2" <0pt>
\ar @{<-} "a";"b_3" <0pt>
\ar @{<-} "a";"b_1" <0pt>
\ar @{<-} "a";"b_4" <0pt>
\ar @{<-} "a";"b_5" <0pt>
\ar @{->} "c1";"b_3" <0pt>
\ar @{->} "c2";"b_3" <0pt>
\ar @{->} "c3";"b_3" <0pt>
\endxy}\Ea
+
\sum_{A\subsetneq [m]\atop \# A\geq 2} 
\Ba{c}\resizebox{22mm}{!}{\xy
(1.6,7)*{...},
(-13,7)*{},
(-8,7)*{},
(-3,7)*{},
(7,7)*{},
(13,7)*{},
 (0,0)*{\bu}="a",
 (0,-6)*{...},
(-5,-7)*{}="01",
(-3,-7)*{}="02",
(3,-7)*{}="03",
(5,-7)*{}="04",
(-12,7)*{}="b_1",
(-9,7)*{}="b_2",
(-6,7)*{...},
(-3,7)*{\bu}="b_3",
(-9,14)*{}="c1",
(-5,14)*{}="c2",
(-1,14)*{...},
(3,14)*{}="c3",
(7,7)*{}="b_4",
(11,7)*{}="b_5",
(-3,17)*{\overbrace{\ \ \ \ \ \ \ \ \ \ \ \ }^A},
\ar @{<-} "a";"01" <0pt>
\ar @{<-} "a";"02" <0pt>
\ar @{<-} "a";"03" <0pt>
\ar @{<-} "a";"04" <0pt>
\ar @{->} "a";"b_2" <0pt>
\ar @{->} "a";"b_3" <0pt>
\ar @{->} "a";"b_1" <0pt>
\ar @{->} "a";"b_4" <0pt>
\ar @{->} "a";"b_5" <0pt>
\ar @{<-} "c1";"b_3" <0pt>
\ar @{<-} "c2";"b_3" <0pt>
\ar @{<-} "c3";"b_3" <0pt>
\endxy}\Ea
+ \ \  \ldots.
\Eeq
We show above explicitly only those terms in the differential which have
the maximal possible number of directed paths from the 
in-legs to the out-legs (that number equals to $mn$; all the omitted terms have this number $<mn$).
 We use in our proof only these terms. 

\sip

The associated cohomology properad $\cI\cB:=H^\bu(\cI\cB_\infty)$ was 
proven in \cite{A} to be generated by the following planar trivalent corollas of degree +1,
\Beq\label{3: generators of IB}
 \Ba{c}\resizebox{8mm}{!}{  \xy
(0,6)*{}="1";
    (0,0.2)*{\bu}="L";
  (-4,-5)*{}="C";
   (+4,-5)*{}="D";
\ar @{->} "D";"L" <0pt>
\ar @{->} "C";"L" <0pt>
\ar @{<-} "1";"L" <0pt>
 \endxy}
 \Ea,
 \ \ \
 \Ba{c}\resizebox{8mm}{!}{  \xy
(0,-6)*{}="1";
    (0,-0.2)*{\bu}="L";
  (-4,5)*{}="C";
   (4,5)*{}="D";
\ar @{<-} "D";"L" <0pt>
\ar @{<-} "C";"L" <0pt>
\ar @{->} "1";"L" <0pt>
 \endxy}
 \Ea, 
\Eeq
modulo the following relations
$$
  \Ba{c}\resizebox{6.5mm}{!}{  \xy
(-4,10)*{}="1";
 (4,10)*{}="2";
    (0,3.5)*{\bbu}="A";
 (0,-3.5)*{\bbu}="B";
 (-4,-10)*{}="b1";
 (4,-10)*{}="b2";
\ar @{->} "A";"1" <0pt>
\ar @{->} "A";"2" <0pt>
\ar @{->} "B";"A" <0pt>
\ar @{<-} "B";"b1" <0pt>
\ar @{<-} "B";"b2" <0pt>
 \endxy}
 \Ea
 +
 \Ba{c}\resizebox{13.5mm}{!}{  \xy
(-9,10)*{}="1";
    (-9,+3)*{\bbu}="L";
 (-14,-3.5)*{\bbu}="B";
 (-19,5)*+{}="b1";
 (-14,-12)*+{}="b2";
  (-3,-5)*{}="C";
\ar @{->} "C";"L" <0pt>
\ar @{->} "B";"L" <0pt>
\ar @{->} "B";"b1" <0pt>
\ar @{<-} "B";"b2" <0pt>
\ar @{<-} "1";"L" <0pt>
 \endxy}
 \Ea
 +
 \Ba{c}\resizebox{13.5mm}{!}{  \xy
(-9,-10)*{}="1";
    (-9,-3)*{\bbu}="L";
 (-14,3.5)*{\bbu}="B";
 (-19,-5)*+{}="b1";
 (-14,12)*+{}="b2";
  (-3,5)*{}="C";
\ar @{<-} "C";"L" <0pt>
\ar @{<-} "B";"L" <0pt>
\ar @{<-} "B";"b1" <0pt>
\ar @{->} "B";"b2" <0pt>
\ar @{->} "1";"L" <0pt>
 \endxy}
 \Ea
=0,
\Ba{c}\resizebox{13.5mm}{!}{  \xy
(-9,10)*{^0}="1";
    (-9,+3)*{\bbu}="L";
 (-14,-3.5)*{\bbu}="B";
 (-20,-12)*+{_1}="b1";
 (-8,-12)*+{_2}="b2";
  (-3,-4)*{_3}="C";
\ar @{->} "C";"L" <0pt>
\ar @{->} "B";"L" <0pt>
\ar @{<-} "B";"b1" <0pt>
\ar @{<-} "B";"b2" <0pt>
\ar @{<-} "1";"L" <0pt>
 \endxy}
 \Ea
 +
  \Ba{c}\resizebox{13.5mm}{!}{  \xy
(9,10)*{^0}="1";
    (9,+3)*{\bbu}="L";
 (14,-3.5)*{\bbu}="B";
 (20,-12)*+{_3}="b1";
 (8,-12)*+{2}="b2";
  (3,-4)*{1}="C";
\ar @{->} "C";"L" <0pt>
\ar @{->} "B";"L" <0pt>
\ar @{<-} "B";"b1" <0pt>
\ar @{<-} "B";"b2" <0pt>
\ar @{<-} "1";"L" <0pt>
 \endxy}
 \Ea=0, 
 \Ba{c}\resizebox{13.5mm}{!}{  \xy
(-9,-10)*{^0}="1";
    (-9,-3)*{\bbu}="L";
 (-14,3.5)*{\bbu}="B";
 (-20,12)*+{_1}="b1";
 (-8,12)*+{_2}="b2";
  (-3,4)*{_3}="C";
\ar @{->} "C";"L" <0pt>
\ar @{->} "B";"L" <0pt>
\ar @{<-} "B";"b1" <0pt>
\ar @{<-} "B";"b2" <0pt>
\ar @{<-} "1";"L" <0pt>
 \endxy}
 \Ea
 +
  \Ba{c}\resizebox{13.5mm}{!}{  \xy
(9,-10)*{^0}="1";
    (9,-3)*{\bbu}="L";
 (14,3.5)*{\bbu}="B";
 (20,12)*+{_3}="b1";
 (8,12)*+{2}="b2";
  (3,4)*{1}="C";
\ar @{<-} "C";"L" <0pt>
\ar @{<-} "B";"L" <0pt>
\ar @{->} "B";"b1" <0pt>
\ar @{->} "B";"b2" <0pt>
\ar @{->} "1";"L" <0pt>
 \endxy}
 \Ea=0.
$$
Strictly speaking, what we discussed just above is a degree shifted (``odd") version of the dg properad introduced and studied in \cite{A}.

\sip

For future reference we note that if we consider a filtration of $\cI\cB_\infty$ by the number of directed paths from in-legs to out-legs, then the cohomology of the associated graded complex $gr \cI\cB_\infty$ is given by a properad $gr \cI\cB$ with the same generators (\ref{3: generators of IB}) but with the first relation in the above list simplified to
$$
  \Ba{c}\resizebox{6.5mm}{!}{  \xy
(-4,10)*{}="1";
 (4,10)*{}="2";
    (0,3.5)*{\bbu}="A";
 (0,-3.5)*{\bbu}="B";
 (-4,-10)*{}="b1";
 (4,-10)*{}="b2";
\ar @{->} "A";"1" <0pt>
\ar @{->} "A";"2" <0pt>
\ar @{->} "B";"A" <0pt>
\ar @{<-} "B";"b1" <0pt>
\ar @{<-} "B";"b2" <0pt>
 \endxy}
 \Ea=0.
$$

\subsection{On vertices of ribbon quivers of valency $\geq 4$} Consider a generic
vertex of valency $\geq 4$  of a ribbon quiver $\Ga\in \ORGC_1$, say this one,
$$
v=  \Ba{c}\resizebox{13mm}{!}{  \xy
(0,5)*{}="U";
(-2,4)*{}="UL";
(2,4)*{}="UR";
(0,-5)*{}="D";
    (0,0)*{\bu}="C";
(-5,1)*{}="L1";
(-5,-1)*{}="L2";    
  (5,0)*{}="R";
  
\ar @{->} "C";"U" <0pt>
\ar @{->} "C";"UL" <0pt>
\ar @{->} "C";"UR" <0pt>
\ar @{->} "C";"D" <0pt>
\ar @{<-} "C";"L1" <0pt>
\ar @{<-} "C";"L2" <0pt>
\ar @{<-} "C";"R" <0pt>
 \endxy}
 \Ea.
$$
All edges attached to $v$ are cyclically ordered. It may happen that two neighboring (with respect to the given cyclic order) edges have the same direction; we call such edges {\it parallel at $v$}. For example, the above vertex has three parallel out-edges and two parallel in-edges. Parallel edges come in maximal cyclically ordered blocks which we call {\it bunches}\, at $v$. A vertex $v$ with no parallel edges has even valency as the total number $n\geq 2$ of incoming edges must be equal to the number of outgoing edges. Here are examples of ribbon vertices  
with {\it no}\, parallel edges,
\Beq\label{3: DLie corollas}
  \Ba{c}\resizebox{13mm}{!}{  \xy
(0,-5)*{}="U";
(0,+5)*{}="D";
    (0,0)*{\bu}="C";
  (5,0)*{}="R";
   (-5,0)*{}="L";
\ar @{->} "C";"U" <0pt>
\ar @{->} "C";"D" <0pt>
\ar @{<-} "C";"L" <0pt>
\ar @{<-} "C";"R" <0pt>
 \endxy}
 \Ea
 ,\
   \Ba{c}\resizebox{13mm}{!}{  \xy
   (-3,4)*{}="UL";
(3,4)*{}="UR";
(-3,-4)*{}="DL";
(3,-4)*{}="DR";
    (0,0)*{\bu}="C";
  (5,0)*{}="R";
   (-5,0)*{}="L";
\ar @{->} "C";"UL" <0pt>
\ar @{<-} "C";"UR" <0pt>
\ar @{->} "C";"DL" <0pt>
\ar @{<-} "C";"DR" <0pt>
\ar @{<-} "C";"L" <0pt>
\ar @{->} "C";"R" <0pt>
 \endxy}
 \Ea.
\Eeq

\subsection{Reminder on the dg properad $\caD\caL ie_\infty$} Let $\caD\caL ie_\infty$ be a graded vector space spanned by connected ribbon quivers $\Ga$  having vertices of valency 1 and of valencies  $\geq 4$ {\it with no parallel edges}. The vertices of valency 1 are of two types, depending on the direction of the attached edge. The vertices of valency 1 (resp.\ $\geq 4$) are assigned the cohomological degree $0$ (resp.\ $+1$).  We call from now one 1-valent vertices {\it hairs}\, (and show them in pictures as hairs or free legs); we also assume that 1-valent vertices are labelled. As vertices of $\Ga$ with valencies $\geq 4$ have no parallel edges or hairs attached, the number, say $n$, of in-hairs of $\Ga$ must be equal to the number of out-hairs so that the permutation group $\bS_n \times \bS_n$ acts on such haired ribbon quivers. The subspace of  $\caD\caL ie_\infty$ spanned by ribbon quivers with $n$ in-hairs (and hence $n$-out hairs) is denoted by $\caD\caL ie_\infty(n,n)$. Clearly, the collection
$$
\caD\caL ie_\infty=\{\caD\caL ie_\infty(n,n)\}_{n\geq 2} 
$$  
is a free properad generated by at least 4-valent {\it ribbon}\, corollas with no parallel edges as in (\ref{3: DLie corollas}), for example 
$$
  \Ba{c}\resizebox{20mm}{!}{  \xy
(0,7)*+{_{\bar{1}}}="U";
(0,-7)*+{_{\bar{2}}}="D";
    (0,0)*{\bu}="C";
  (7,0)*+{_2}="R";
   (-7,0)*+{_1}="L";
\ar @{->} "C";"U" <0pt>
\ar @{->} "C";"D" <0pt>
\ar @{<-} "C";"L" <0pt>
\ar @{<-} "C";"R" <0pt>
 \endxy}
 \Ea\in  \caD\caL ie_\infty(2,2) \ , \ 
 \Ba{c}\resizebox{28mm}{!}{  \xy
(0,-7)*+{_{\bar{3}}}="D";
(0,+7)*+{_{\bar{1}}}="U";
 (-7,0)*+{_1}="L";
(0,0)*{\bu}="C";
(7,0)*{\bu}="CR";
  (7,-7)*+{_3}="RD";
(7,+7)*+{_2}="RU";
 (14,0)*+{_{\bar{2}}}="RR";
\ar @{->} "C";"U" <0pt>
\ar @{->} "C";"D" <0pt>
\ar @{<-} "C";"L" <0pt>
\ar @{<-} "C";"CR" <0pt>
\ar @{<-} "CR";"RU" <0pt>
\ar @{<-} "CR";"RD" <0pt>
\ar @{->} "CR";"RR" <0pt>
 \endxy}
 \Ea \in \caD\caL ie_\infty(3,3)
$$
The differential $\delta$ in $\mathcal{DL}ie_\infty$ acts on such an $(n,n)$-corolla by substituting into the vertex the graph $\Ba{c}\xy
 (0,0)*{\bu}="a",
(0,5)*{\bu}="b",
\ar @{->} "a";"b" <0pt>
\endxy\Ea$ and reattaching the edges among the two newly created vertices in such a way that the cyclic structure is preserved and the newly created vertices are at least 4-valent and have no parallel edges attached, e.g.
\Beq\label{3: d in DLie_infty example}
\delta \Ba{c}\resizebox{16mm}{!}{  \xy
   (-5,6)*{_1}="UL";
(5,6)*{_2}="UR";
(-5,-6)*{_5}="DL";
(5,-6)*{_4}="DR";
    (0,0)*{\bu}="C";
  (7,0)*{_3}="R";
   (-7,0)*{_0}="L";
\ar @{->} "C";"UL" <0pt>
\ar @{<-} "C";"UR" <0pt>
\ar @{->} "C";"DL" <0pt>
\ar @{<-} "C";"DR" <0pt>
\ar @{<-} "C";"L" <0pt>
\ar @{->} "C";"R" <0pt>
 \endxy}
 \Ea
 =
   \Ba{c}\resizebox{28mm}{!}{  \xy
(0,-7)*+{_5}="D";
(0,+7)*+{_1}="U";
 (-7,0)*+{_0}="L";
(0,0)*{\bu}="C";
(7,0)*{\bu}="CR";
  (7,-7)*+{_4}="RD";
(7,+7)*+{_2}="RU";
 (14,0)*+{_3}="RR";
\ar @{->} "C";"U" <0pt>
\ar @{->} "C";"D" <0pt>
\ar @{<-} "C";"L" <0pt>
\ar @{<-} "C";"CR" <0pt>
\ar @{<-} "CR";"RU" <0pt>
\ar @{<-} "CR";"RD" <0pt>
\ar @{->} "CR";"RR" <0pt>
 \endxy}
 \Ea
 +
  \Ba{c}\resizebox{28mm}{!}{  \xy
(0,-7)*+{_3}="D";
(0,+7)*+{_5}="U";
 (-7,0)*+{_4}="L";
(0,0)*{\bu}="C";
(7,0)*{\bu}="CR";
  (7,-7)*+{_2}="RD";
(7,+7)*+{_0}="RU";
 (14,0)*+{_1}="RR";
\ar @{->} "C";"U" <0pt>
\ar @{->} "C";"D" <0pt>
\ar @{<-} "C";"L" <0pt>
\ar @{<-} "C";"CR" <0pt>
\ar @{<-} "CR";"RU" <0pt>
\ar @{<-} "CR";"RD" <0pt>
\ar @{->} "CR";"RR" <0pt>
 \endxy}
 \Ea
  +
  \Ba{c}\resizebox{28mm}{!}{  \xy
(0,-7)*+{_1}="D";
(0,+7)*+{_3}="U";
 (-7,0)*+{_2}="L";
(0,0)*{\bu}="C";
(7,0)*{\bu}="CR";
  (7,-7)*+{_0}="RD";
(7,+7)*+{_4}="RU";
 (14,0)*+{_5}="RR";
\ar @{->} "C";"U" <0pt>
\ar @{->} "C";"D" <0pt>
\ar @{<-} "C";"L" <0pt>
\ar @{<-} "C";"CR" <0pt>
\ar @{<-} "CR";"RU" <0pt>
\ar @{<-} "CR";"RD" <0pt>
\ar @{->} "CR";"RR" <0pt>
 \endxy}
 \Ea.
\Eeq

It was proven in \cite{L} that the cohomology properad $\caD\caL ie:=H^\bu(\caD\caL ie)$
is generated by 4-valent ribbon corolla   $\Ba{c}\resizebox{11mm}{!}{  \xy
(0,7)*{}="U";
(0,-7)*{}="D";
    (0,0)*{\bbu}="C";
  (7,0)*{}="R";
   (-7,0)*{}="L";
\ar @{->} "C";"U" <0pt>
\ar @{->} "C";"D" <0pt>
\ar @{<-} "C";"L" <0pt>
\ar @{<-} "C";"R" <0pt>
 \endxy}
 \Ea$ modulo the relations (\ref{2: double Lie relations}). Representations of the properad $\caD\caL ie$ are (degree shifted) double Lie algebras introduced and studied in \cite{vdB} as a part of a richer double Poisson structure.

 \subsubsection{\bf Remark} Strictly speaking, the properad $\caD\caL ie$   we discussed just above is a degree shifted (``odd") version of the properad $\caD\caL ie^{even}$ studied in \cite{L}; the latter properad is concentrated in degree zero. Moreover, the elements of that double Lie properad  $\caD\caL ie^{even}$ and of its minimal resolution $\caD\caL ie^{even}_\infty$ were {not} understood in \cite{L} as {\it ribbon}\, quivers but as elements of a protoperad.
 
 \sip
 
 The original construction in \cite{L} can be understood as generated by $\geq 4$ {\it ribbon}\, corollas with no parallel edges whose vertices are all assigned the cohomological degree +2, the outgoing half-edges (or out-hairs) are assigned the degree $-1$ and the in-going half-edges (or in-hairs) are assigned the degree $0$ (so that the internal edges of the elements of $\caD\caL ie^{even}_\infty$  all have the overall degree $-1$
 corresponding to the value $d=2$ of our integer parameter). In this even case the orientation of a ribbon quiver generator of  $\caD\caL ie^{even}_\infty$  is given by an ordering of its edges and out-hairs (rather than by ordering of the vertices as in our odd version discussed above). By the way, in this  ribbon quivers incarnation of the protoperad of even double Lie algebras  $\caD\caL ie^{even}$  the skew-symmetry of the generating 4-valent  ribbon corolla  becomes automatic, it follows from cyclic ordering of attached hairs and from the definition of its orientation, i.e. there is no need to impose this skew-symmetry as an extra condition.

  \subsection{Remark on convergence of spectral sequences} We study below several complexes $C$ which are either equal to $\ORGC_1$ or are its quotients by homogeneous (with respect to the genus and the number of boundaries) subcomplexes. Every such a complex $C$ decomposes into a direct sum of complexes 
  $$
  C= \prod_{g\geq 0, m\geq 1\atop 2g+m\geq 3} C^{g,m}
  $$
parameterized by the genus $g$ and the number of boundaries $m$ of the generating (equivalence classes of)  ribbon quivers;
this decomposition in all cases is inherited from $\ORGC_1$.
The cohomological degree of a generator $\Ga$ in each complex $C^{g,m}$ is given by one and the same formula 
$$
|\Ga|=\# V(\Ga)=\text{the number of vertices of $\Ga$}.
$$
We also have a relation
$$
\#E(\Ga) - \# V(\Ga)= 2g-2+m.
$$
Thus if $|\Ga|$ is fixed, then $\# V(\Ga)$ and $\#E(\Ga)$ are also fixed. We conclude that a linear subspace of $C^{g,m}$ spanned by ribbon quivers with a fixed cohomological degree is finite-dimensional. Hence any filtration of $C^{g,m}$  which respects the differential is bounded, and hence it converges. We consider below only those filtrations of $C$ which induce well-defined filtrations on each summand $C^{g,m}$; the associated spectral sequences  all converge.

\subsection{First reduction of $\ORGC_{1}$} Let $K$ be a linear subspace 
of the complex $\ORGC_{1}$ spanned by ribbon quivers having at least one trivalent source or trivalent target, or at least one vertex $v$ of valency $\geq 4$ which has at least one pair of parallel edges at $v$, and let $\langle K, \delta K\rangle$ be its differential closure in $\ORGC_{1}$. Define a ribbon graph complex  $\widehat{\Delta}\RGC_{1}$  by the following short exact sequence
$$
0 \lon \langle K, \delta K\rangle \lon \ORGC_{1} \stackrel{p}{\lon} \widehat{\Delta}\RGC_{1} \lon 0.
$$
The quotient complex  $\widehat{\Delta}\RGC_{1}$ is generated by equivalence classes of ribbon quivers which can have 2- and 3-valent vertices of the form
\Beq\label{2: 2 and 3 valent vertices}
\Ba{c}\resizebox{8mm}{!}{ \xy
 (0,0)*{}="a",
(4,3)*{\bullet}="b",
(8,0)*{}="c",
\ar @{<-} "a";"b" <0pt>
\ar @{->} "b";"c" <0pt>
\endxy}\Ea, 
\ \ \ \  
\Ba{c}\resizebox{8mm}{!}{\xy
 (0,0)*{}="a",
(4,3)*{\bullet}="b",
(8,0)*{}="c",
\ar @{->} "a";"b" <0pt>
\ar @{<-} "b";"c" <0pt>
\endxy} \Ea, 
\ \ \ 
 \Ba{c}\resizebox{8mm}{!}{  \xy
(0,6)*{}="1";
    (0,0.2)*{\bu}="L";
  (-4,-5)*{}="C";
   (+4,-5)*{}="D";
\ar @{->} "D";"L" <0pt>
\ar @{->} "C";"L" <0pt>
\ar @{<-} "1";"L" <0pt>
 \endxy}
 \Ea,
 \ \ \
 \Ba{c}\resizebox{8mm}{!}{  \xy
(0,-6)*{}="1";
    (0,-0.2)*{\bu}="L";
  (-4,5)*{}="C";
   (4,5)*{}="D";
\ar @{<-} "D";"L" <0pt>
\ar @{<-} "C";"L" <0pt>
\ar @{->} "1";"L" <0pt>
 \endxy}
 \Ea, 
\Eeq
and also vertices of valencies $2k$, $k\geq 2$, with no parallel edges attached (as, e.g., in (\ref{3: DLie corollas})).
The latter set of vertices can be identified with the set of ribbon generators of the dg properad 
$\caD\Lie_\infty$. For each internal edge of a ribbon quiver $\Ga\in  \widehat{\Delta}\RGC_{1}$  there are relations  (\ref{2: 2 and 3 valent vertices}), 
(\ref{2: twisted ass-type relations}), (\ref{2: IB type relation}) and also relations of the type  (\ref{2: 3+4 relations})  for each edge connecting a trivalent vertex to a $2k$-valent vertex of $\caD\Lie_\infty$-type in the full analogy to the case  $k=2$.

\sip

The induced differential $\delta$ in $\widehat{\Delta}\RGC_{1}$ acts on vertices of valencies $\leq 4$ as in 
(\ref{2 d on 2val})-(\ref{2: d on 4val}), and $\delta$ acts on $2k$-valent vertices with $k\geq 3$ 
precisely as the diffenetial in $\caD\Lie_\infty$ acts on its ribbon generators, see (\ref{3: d in DLie_infty example}) for an example.

\subsubsection{\bf Proposition}\label{3: Proposition on p}  {\it The epimorphism of complexes 
\Beq\label{3: p'}
 \ORGC_{1} \stackrel{p}{\lon} \widehat{\Delta}\RGC_{1} 
\Eeq
 is a quasi isomorphism}.
 
\begin{proof} Consider first filtrations of both sides of the epimorphism $p$
by the number of sources and targets. If we show that the induced map
of the associated graded complexes
$$
gr(p): gr\ORGC_1 \lon  gr \widehat{\Delta}\RGC_{1}
$$
is a quasi-isomorphism, we are done.

\sip

Consider next a filtration of both sides of the morphism $gr(p)$ by
the total number of directed continuous  paths of edges which connect the sources to targets
(the differentials of both sides of $gr(p)$ can not increase this number,  so this filtration is well defined). We obtain a sequence of morphisms
$$
gr(p)_r: (\cE_r gr\ORGC_{1},\delta_r)  \lon (\cE_r gr \widehat{\Delta}\RGC_{1}, d_r), \ \ \ r\in \N,
$$
of the $r$-th pages of the associated spectral sequences. If we show that
$gr(p)_1$ is an isomorphism of complexes, we are done.

\sip

The complex $\cE_0 gr \widehat{\Delta}\RGC_{1}$ is trivial, $d_0=0$, so that 
$\cE_1 gr \widehat{\Delta}\RGC_{1} \simeq \cE_0 gr \widehat{\Delta}\RGC_{1}$ as  graded vector spaces.  It is important to note that the relations on internal edges of generators of $\cE_0  gr \widehat{\Delta}\RGC_{1}$ take a much simpler form 
\Beq\label{3: simplified relations}
  \Ba{c}\resizebox{5.0mm}{!}{  \xy
(-4,10)*{}="1";
 (4,10)*{}="2";
    (0,3.5)*{\bbu}="A";
 (0,-3.5)*{\bbu}="B";
 (-4,-10)*{}="b1";
 (4,-10)*{}="b2";
\ar @{->} "A";"1" <0pt>
\ar @{->} "A";"2" <0pt>
\ar @{->} "B";"A" <0pt>
\ar @{<-} "B";"b1" <0pt>
\ar @{<-} "B";"b2" <0pt>
 \endxy}
 \Ea
\hspace{-2mm} =0,
\Ba{c}\resizebox{10.5mm}{!}{  \xy
(-9,10)*{^0}="1";
    (-9,+3)*{\bbu}="L";
 (-14,-3.5)*{\bbu}="B";
 (-20,-12)*+{_1}="b1";
 (-8,-12)*+{_2}="b2";
  (-3,-4)*{_3}="C";
\ar @{->} "C";"L" <0pt>
\ar @{->} "B";"L" <0pt>
\ar @{<-} "B";"b1" <0pt>
\ar @{<-} "B";"b2" <0pt>
\ar @{<-} "1";"L" <0pt>
 \endxy}
 \Ea
 +
  \Ba{c}\resizebox{10.5mm}{!}{  \xy
(9,10)*{^0}="1";
    (9,+3)*{\bbu}="L";
 (14,-3.5)*{\bbu}="B";
 (20,-12)*+{_3}="b1";
 (8,-12)*+{2}="b2";
  (3,-4)*{1}="C";
\ar @{->} "C";"L" <0pt>
\ar @{->} "B";"L" <0pt>
\ar @{<-} "B";"b1" <0pt>
\ar @{<-} "B";"b2" <0pt>
\ar @{<-} "1";"L" <0pt>
 \endxy}
 \Ea
 \hspace{-2mm}=0, 
 \Ba{c}\resizebox{10.5mm}{!}{  \xy
(-9,-10)*{^0}="1";
    (-9,-3)*{\bbu}="L";
 (-14,3.5)*{\bbu}="B";
 (-20,12)*+{_1}="b1";
 (-8,12)*+{_2}="b2";
  (-3,4)*{_3}="C";
\ar @{->} "C";"L" <0pt>
\ar @{->} "B";"L" <0pt>
\ar @{<-} "B";"b1" <0pt>
\ar @{<-} "B";"b2" <0pt>
\ar @{<-} "1";"L" <0pt>
 \endxy}
 \Ea
 +
  \Ba{c}\resizebox{10.5mm}{!}{  \xy
(9,-10)*{^0}="1";
    (9,-3)*{\bbu}="L";
 (14,3.5)*{\bbu}="B";
 (20,12)*+{_3}="b1";
 (8,12)*+{2}="b2";
  (3,4)*{1}="C";
\ar @{<-} "C";"L" <0pt>
\ar @{<-} "B";"L" <0pt>
\ar @{->} "B";"b1" <0pt>
\ar @{->} "B";"b2" <0pt>
\ar @{->} "1";"L" <0pt>
 \endxy}
 \Ea
 \hspace{-2mm}=0,
   \Ba{c}\resizebox{16mm}{!}{  \xy
(-3,13)*+{_1}="UL";
(3,13)*+{_2}="UR";
(0,7)*{\bu}="U";
(0,-7)*+{_{2k}}="D";
(0,-10)*{\underbrace{\ \ \ \ \ \ \ \ \ \ \ \ \ \ \ \ \ }_{2k\ \text{valent vertex}}};
    (0,0)*{\bu}="C";
  (7,0)*+{_3}="R";
   (-7,0)*+{_0}="L";
\ar @{->} "U";"UL" <0pt>   
\ar @{->} "U";"UR" <0pt>  
\ar @{->} "C";"U" <0pt>
\ar @{->} "C";"D" <0pt>
\ar @{<-} "C";"L" <0pt>
\ar @{<-} "C";"R" <0pt>
 \endxy}
 \Ea
\hspace{-2mm}=0,
   \Ba{c}\resizebox{16mm}{!}{  \xy
(-3,13)*+{_1}="UL";
(3,13)*+{_2}="UR";
(0,7)*{\bu}="U";
(0,-7)*+{_{2k}}="D";
(0,-10)*{\underbrace{\ \ \ \ \ \ \ \ \ \ \ \ \ \ \ \ \ }_{2k\ \text{valent vertex}}};
    (0,0)*{\bu}="C";
  (7,0)*+{_3}="R";
   (-7,0)*+{_0}="L";
\ar @{<-} "U";"UL" <0pt>   
\ar @{<-} "U";"UR" <0pt>  
\ar @{<-} "C";"U" <0pt>
\ar @{<-} "C";"D" <0pt>
\ar @{->} "C";"L" <0pt>
\ar @{->} "C";"R" <0pt>
 \endxy}
 \Ea\hspace{-2mm}=0,
\Eeq
comparing to those in  $\widehat{\Delta}\RGC_{1}$; they all can be easily resolved in terms of a suitable basis in  $\cE_0 gr \widehat{\Delta}\RGC_{1}$.

\sip

The differential $\delta_0$ on the complex  $\cE_0 gr\ORGC_{1}$ acts non-trivially only on (i) trivalent sources and targets and (ii) on those vertices of valency $\geq 4$ which have at least one pair of parallel edges. The edges attached to any vertex $v$ of a generator $\Ga\in \cE_0 gr\ORGC_{1}$ can be grouped in bunches consisting of {\it all}\, neighboring edges parallel to each other; the cardinality of such a bunch is $\geq 1$, and the total number $b(v)$ of such bunches at $v$  is equal to $1$ if $v$ is a source or a target in $\Ga$, or it is equal to an even number otherwise.
 For example, the first two vertices in the following set
$$
v_1=  \Ba{c}\resizebox{13mm}{!}{  \xy
(0,-5)*{}="U";
(0,+5)*{}="D";
    (0,0)*{\bu}="C";
  (5,0)*{}="R";
   (-5,0)*{}="L";
\ar @{->} "C";"U" <0pt>
\ar @{->} "C";"D" <0pt>
\ar @{<-} "C";"L" <0pt>
\ar @{<-} "C";"R" <0pt>
 \endxy}
 \Ea, \ \ 
v_2=  \Ba{c}\resizebox{13mm}{!}{  \xy
(0,5)*{}="U";
(-2,4)*{}="UL";
(2,4)*{}="UR";
(0,-5)*{}="D";
    (0,0)*{\bu}="C";
(-5,1)*{}="L1";
(-5,-1)*{}="L2";    
  (5,0)*{}="R";
  
\ar @{->} "C";"U" <0pt>
\ar @{->} "C";"UL" <0pt>
\ar @{->} "C";"UR" <0pt>
\ar @{->} "C";"D" <0pt>
\ar @{<-} "C";"L1" <0pt>
\ar @{<-} "C";"L2" <0pt>
\ar @{<-} "C";"R" <0pt>
 \endxy}
 \Ea
 ,
 \ \ 
 v_3=
 \Ba{c}\resizebox{11mm}{!}{ \xy
(0,8)*{\overbrace{\  \ \ \ \ \ \ \ \ \ }^{m\geq 2}},
(0,-8)*{\underbrace{\  \ \ \ \ \ \ \ \ \ }_{n\geq 2}},
(0,4.5)*+{...},
(0,-4.5)*+{...},
(0,0)*{\bu}="o",
(-5,5)*{}="1",
(-3,5)*{}="2",
(3,5)*{}="3",
(5,5)*{}="4",
(-3,-5)*{}="5",
(3,-5)*{}="6",
(5,-5)*{}="7",
(-5,-5)*{}="8",
\ar @{->} "o";"1" <0pt>
\ar @{->} "o";"2" <0pt>
\ar @{->} "o";"3" <0pt>
\ar @{->} "o";"4" <0pt>
\ar @{<-} "o";"5" <0pt>
\ar @{<-} "o";"6" <0pt>
\ar @{<-} "o";"7" <0pt>
\ar @{<-} "o";"8" <0pt>
\endxy}\Ea
$$
have the same number of bunches $b(v_1)=b(v_2)=4$, while $b(v_3)=2$.
 Let us call vertices $v$ of $\Ga$ {\it operadic}\, if
$b(v)=2$ and precisely one bunch at $v$ has cardinality 1. Thus operadic vertices have one of the following two forms
$$
  \underbrace{\Ba{c}\resizebox{20mm}{!}
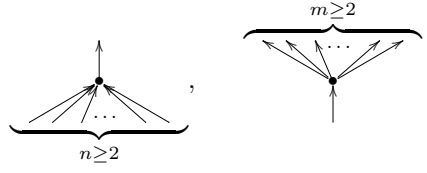
{\xy
(1,-6)*{\ldots},
(-13,-7)*{},
(-8,-7)*{},
(-3,-7)*{},
(7,-7)*{},
(13,-7)*{},
 (0,0)*{\bu}="a",
(0,7)*{}="0",
(-12,-7)*{}="b_1",
(-8,-7)*{}="b_2",
(-3,-7)*{}="b_3",
(8,-7)*{}="b_4",
(12,-7)*{}="b_5",
\ar @{->} "a";"0" <0pt>
\ar @{<-} "a";"b_2" <0pt>
\ar @{<-} "a";"b_3" <0pt>
\ar @{<-} "a";"b_1" <0pt>
\ar @{<-} "a";"b_4" <0pt>
\ar @{<-} "a";"b_5" <0pt>
\endxy}\Ea}_{n\geq 2}, \ \ \ \ \ 
  \overbrace{\Ba{c}\resizebox{20mm}{!}{\xy
(1,6)*{\ldots},
(-13,7)*{},
(-8,7)*{},
(-3,7)*{},
(7,7)*{},
(13,7)*{},
 (0,0)*{\bu}="a",
(0,-7)*{}="0",
(-12,7)*{}="b_1",
(-8,7)*{}="b_2",
(-3,7)*{}="b_3",
(8,7)*{}="b_4",
(12,7)*{}="b_5",
\ar @{<-} "a";"0" <0pt>
\ar @{->} "a";"b_2" <0pt>
\ar @{->} "a";"b_3" <0pt>
\ar @{->} "a";"b_1" <0pt>
\ar @{->} "a";"b_4" <0pt>
\ar @{->} "a";"b_5" <0pt>
\endxy}\Ea}^{m\geq 2}.
$$
Operadic vertices have a distinguished edge which belongs to the bunch
of cardinality one; let us call such an edge {\it operadic}\, as well.

\sip

The main point is that the preserving the path filtration differential $\delta_0$ can create {\it only operadic vertices}, for example
$$
\delta_0
\Ba{c}\resizebox{11mm}{!}{ \xy
(0,8)*{\overbrace{\  \ \ \ \ \ \ \ \ \ }^{m\geq 2}},
(0,-8)*{\underbrace{\  \ \ \ \ \ \ \ \ \ }_{n\geq 2}},
(0,4.5)*+{...},
(0,-4.5)*+{...},
(0,0)*{\bu}="o",
(-5,5)*{}="1",
(-3,5)*{}="2",
(3,5)*{}="3",
(5,5)*{}="4",
(-3,-5)*{}="5",
(3,-5)*{}="6",
(5,-5)*{}="7",
(-5,-5)*{}="8",
\ar @{->} "o";"1" <0pt>
\ar @{->} "o";"2" <0pt>
\ar @{->} "o";"3" <0pt>
\ar @{->} "o";"4" <0pt>
\ar @{<-} "o";"5" <0pt>
\ar @{<-} "o";"6" <0pt>
\ar @{<-} "o";"7" <0pt>
\ar @{<-} "o";"8" <0pt>
\endxy}\Ea
=
\Ba{c}\resizebox{11mm}{!}{ \xy
(0,9.5)*+{...},
(0,-4.5)*+{...},
(0,5)*{\bu}="oo",
(0,0)*{\bu}="o",
(-5,10)*{}="1",
(-3,10)*{}="2",
(3,10)*{}="3",
(5,10)*{}="4",
(-3,-5)*{}="5",
(3,-5)*{}="6",
(5,-5)*{}="7",
(-5,-5)*{}="8",
\ar @{->} "oo";"1" <0pt>
\ar @{->} "oo";"2" <0pt>
\ar @{->} "oo";"3" <0pt>
\ar @{->} "oo";"4" <0pt>
\ar @{->} "o";"oo" <0pt>
\ar @{<-} "o";"5" <0pt>
\ar @{<-} "o";"6" <0pt>
\ar @{<-} "o";"7" <0pt>
\ar @{<-} "o";"8" <0pt>
\endxy}\Ea
+
\sum_{A\subset [n]\atop \# A\geq 2} 
\Ba{c}\resizebox{22mm}{!}{\xy
(1.6,-7)*{...},
(-13,-7)*{},
(-8,-7)*{},
(-3,-7)*{},
(7,-7)*{},
(13,-7)*{},
 (0,0)*{\bu}="a",
 (0,6)*{...},
(-5,7)*{}="01",
(-3,7)*{}="02",
(3,7)*{}="03",
(5,7)*{}="04",
(-12,-7)*{}="b_1",
(-9,-7)*{}="b_2",
(-6,-7)*{...},
(-3,-7)*{\bu}="b_3",
(-9,-14)*{}="c1",
(-5,-14)*{}="c2",
(-1,-14)*{...},
(3,-14)*{}="c3",
(7,-7)*{}="b_4",
(11,-7)*{}="b_5",
(-3,-17)*{\underbrace{\ \ \ \ \ \ \ \ \ \ \ \ }_A},
\ar @{->} "a";"01" <0pt>
\ar @{->} "a";"02" <0pt>
\ar @{->} "a";"03" <0pt>
\ar @{->} "a";"04" <0pt>
\ar @{<-} "a";"b_2" <0pt>
\ar @{<-} "a";"b_3" <0pt>
\ar @{<-} "a";"b_1" <0pt>
\ar @{<-} "a";"b_4" <0pt>
\ar @{<-} "a";"b_5" <0pt>
\ar @{->} "c1";"b_3" <0pt>
\ar @{->} "c2";"b_3" <0pt>
\ar @{->} "c3";"b_3" <0pt>
\endxy}\Ea
+
\sum_{A\subset [m]\atop \# A\geq 2} 
\Ba{c}\resizebox{22mm}{!}{\xy
(1.6,7)*{...},
(-13,7)*{},
(-8,7)*{},
(-3,7)*{},
(7,7)*{},
(13,7)*{},
 (0,0)*{\bu}="a",
 (0,-6)*{...},
(-5,-7)*{}="01",
(-3,-7)*{}="02",
(3,-7)*{}="03",
(5,-7)*{}="04",
(-12,7)*{}="b_1",
(-9,7)*{}="b_2",
(-6,7)*{...},
(-3,7)*{\bu}="b_3",
(-9,14)*{}="c1",
(-5,14)*{}="c2",
(-1,14)*{...},
(3,14)*{}="c3",
(7,7)*{}="b_4",
(11,7)*{}="b_5",
(-3,17)*{\overbrace{\ \ \ \ \ \ \ \ \ \ \ \ }^A},
\ar @{<-} "a";"01" <0pt>
\ar @{<-} "a";"02" <0pt>
\ar @{<-} "a";"03" <0pt>
\ar @{<-} "a";"04" <0pt>
\ar @{->} "a";"b_2" <0pt>
\ar @{->} "a";"b_3" <0pt>
\ar @{->} "a";"b_1" <0pt>
\ar @{->} "a";"b_4" <0pt>
\ar @{->} "a";"b_5" <0pt>
\ar @{<-} "c1";"b_3" <0pt>
\ar @{<-} "c2";"b_3" <0pt>
\ar @{<-} "c3";"b_3" <0pt>
\endxy}\Ea
$$
 Hence we can argue as follows: given any generator $\Ga\in \cE_0 gr\ORGC_{1}$, let $\widetilde{\Ga}$ be the ribbon quiver obtained from $\Ga$ by contracting its all operadic edges; note that this procedure  does {\it not}\, create {\it closed}\, paths of directed edges, i.e.\ $\widetilde{\Ga}$ is a ribbon {\it quiver}\, indeed. Let us denote the set of such reduced ribbon quivers without operadic vertices by $\widetilde{S}_{red}$. 

\sip

For any $G\in \widetilde{S}_{red}$ denote by  $\cE_0 gr\ORGC_{1}^{G}$
the linear subspace of $\cE_0 gr\ORGC_{1}$ spanned by all ribbon quivers 
$\Ga$ satisfying the condition $\widetilde{\Ga}=G$. As $\widetilde{(\delta_0\Ga)}=\widetilde{\Ga}$, this subspace is a subcomplex; moreover we have a direct sum decomposition
$$
(\cE_0 gr\ORGC_{1}, \delta_0)=  \bigoplus_{G\in \widetilde{S}_{red}} (\cE_0 gr\ORGC_{1}^{G}, \delta_0).
$$
Similarly one one has a direct sum decomposition of graded vector spaces,
$$
\cE_1 gr\widehat{\Delta}\RGC_{1} = \cE_0 gr \widehat{\Delta}\RGC_{1}= \bigoplus_{G\in \widetilde{S}_{red}} \cE_0 gr\widehat{\Delta}\RGC_{1}^{G},
$$ 
The Proposition is proven once we show that for any $G\in \widetilde{S}_{red}$ one has an isomorphism of graded vector spaces 
\Beq\label{3: E_1 G isomorphism}
 H^\bu(\cE_0 gr\ORGC_{1}^{G}, \delta_0)\simeq  \cE_0 gr \widehat{\Delta}\RGC_{1}^{G}
\Eeq
as this implies an isomorphism
$$
\cE_1 gr\ORGC_{1} = \cE_1 gr\widehat{\Delta}\RGC_{1}
$$
(the induced differentials $\delta_1$ and $d_1$ on $\cE_1$-pages must be identical to each other by the very definition of $\widehat{\Delta}\RGC_{1}$).

\sip

For a ribbon quiver $G\in  \widetilde{S}_{red}$ denote by
\Bi
\item[(1)] $V_s(G)$ (resp.\ $V_t(G)$) the set of sources (resp.\ targets) of $G$; the valency of
a vertex $v\in V_s(G)$ (resp., $v\in V_t(G)$) is denoted by $|v|_{out}$ (resp.\ by $|v|_{in}$).

\item[(2)] $V_2(G)$ the set of vertices with $b(v)=2$; such a vertex has $|v|_{in}\geq 2$ incoming edges and  $|v|_{out}\geq 2$  outgoing edges;
\item[(3)] $V_{>2}(G)$ the set of vertices with $b(v)>2$; the set of bunches at $v$ with outgoing (resp., ingoing)  edges is denoted by $B_{out}(v)$ (resp.,  $B_{in}(v)$),
    and the cardinality of a bunch $b$ in   $B_{out}(v)$ (resp., in  $B_{in}(v)$) is denoted by $|b|$.
\Ei
Using this notation we can decompose the complex $(\cE_0 gr\ORGC_{1}^{G}, \delta_0)$
into the unordered\footnote{Using Maschke's theorem, one can assume without loss of generality that the tensor product is the standard (totally ordered) one.} tensor product of the complexes we discussed above,
\Beqrn
\cE_0 gr\ORGC_{1}^{G}& \simeq & \left(\bigotimes_{v\in V_t(G)}  \cA ss_\infty^\wedge(|v|_{in})\right)\ot \left( \bigotimes_{v\in V_s(G)}  \cA ss_\infty^\vee(|v|_{out})\right) \ot \left( \bigotimes_{v\in V_2(G)} gr \cI\cB_\infty (|v|_{out}, |v|_{in})\right)\ot \\
&& \ot \left(\bigotimes_{v\in V_{>2}(G)}\left( \bigotimes_{b\in B_{in}(v)}  \cA ss_\infty^{\wedge,+}(|b|)  \bigotimes_{b\in B_{out}(v)}  \cA ss_\infty^{\vee,+}(|b|)        \right)\right)
\Eeqrn
Using now the well-known results (described in detail in \S\S\  3.2-3.4) on cohomology groups of the above tensor factors, one obtains the isomorphism (\ref{3: E_1 G isomorphism})  implying that the map $gr(p)_1$ of first pages of spectral sequences is an isomorphism o9f complexes. By the classical comparison theorem, The Proposition follows.
\end{proof}

\subsection{Second reduction of $\ORGC_{1}$}
Let $L$ be a linear subspace 
of the complex  $\widehat{\Delta}\RGC_{1}$ spanned by (equivalence classes of) ribbon quivers  having at least one 
vertex $v$ of valency $>4$, and let $\langle L, \delta L\rangle$ be its differential closure in  $\widehat{\Delta}\RGC_{1}$. There is  short exact sequence of complexes
$$
0 \lon \langle L, \delta L\rangle \lon \widehat{\Delta}\RGC_{1}\stackrel{q}{\lon} 
\sRGC_1 \lon 0.
$$

\subsubsection{\bf Proposition} {\it The epimorphism of complexes
$$
 \widehat{\Delta}\RGC_{1}\stackrel{q}{\lon} 
\sRGC_1 
$$ 
is a quasi-isomorphism}.

\begin{proof}
 Consider a filtration of both sides of the morphism $q$ 
by the number of 2- and 3-valent vertices, and let
$$
q_r: \cE_r (\widehat{\Delta}\RGC_{1}) \lon \cE_r (\sRGC_1)
$$
be the induced set of morphisms of $r$th pages of the associated spectral sequences. If we show that
\Beq\label{3: map q_1}
q_1: \cE_1 (\widehat{\Delta}\RGC_{1}) \lon \cE_1 (\sRGC_1)
\Eeq
is an isomorphism of complexes, we are done.

\sip

As the complex $\cE_0 (\sRGC_1)$ has trivial differential, we have an isomorphism
of graded vector spaces,
$$
\cE_0 (\sRGC_1)\simeq \cE_1 (\sRGC_1).
$$
The differential on the complex $\cE_0(\widehat{\Delta}\RGC_{1})$ acts
non-trivially only on vertices of  valencies $>4$ splitting them precisely as the differential in the dg free properad $\caD\Lie_\infty$ splits its generators (see (\ref{3: d in DLie_infty example}) for an example).

\sip

The argument using a path filtration in the proof of the Proposition {\ref{3: Proposition on p}} implies that the equivalence classes of ribbon
quivers generating $\cE_0(\widehat{\Delta}\RGC_{1})$ (or  $\cE_0({\Delta}\RGC_{1})$) admit representatives given by ribbon quivers whose internal edges satisfy the relations (\ref{3: simplified relations}); put another way, we may assume without loss of generality
 that both complexes $\cE_0(\widehat{\Delta}\RGC_{1})$ and  $\cE_0({\Delta}\RGC_{1})$ are spanned by ribbon quivers with {\sip no}\, internal edges of the form
$$
 \Ba{c}\resizebox{5.0mm}{!}
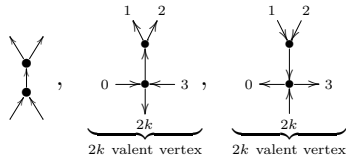
{  \xy
(-4,10)*{}="1";
 (4,10)*{}="2";
    (0,3.5)*{\bbu}="A";
 (0,-3.5)*{\bbu}="B";
 (-4,-10)*{}="b1";
 (4,-10)*{}="b2";
\ar @{->} "A";"1" <0pt>
\ar @{->} "A";"2" <0pt>
\ar @{->} "B";"A" <0pt>
\ar @{<-} "B";"b1" <0pt>
\ar @{<-} "B";"b2" <0pt>
 \endxy}
 \Ea,
    \Ba{c}\resizebox{16mm}{!}{  \xy
(-3,13)*+{_1}="UL";
(3,13)*+{_2}="UR";
(0,7)*{\bu}="U";
(0,-7)*+{_{2k}}="D";
(0,-10)*{\underbrace{\ \ \ \ \ \ \ \ \ \ \ \ \ \ \ \ \ }_{2k\ \text{valent vertex}}};
    (0,0)*{\bu}="C";
  (7,0)*+{_3}="R";
   (-7,0)*+{_0}="L";
\ar @{->} "U";"UL" <0pt>   
\ar @{->} "U";"UR" <0pt>  
\ar @{->} "C";"U" <0pt>
\ar @{->} "C";"D" <0pt>
\ar @{<-} "C";"L" <0pt>
\ar @{<-} "C";"R" <0pt>
 \endxy}
 \Ea
\hspace{-2mm},
   \Ba{c}\resizebox{16mm}{!}{  \xy
(-3,13)*+{_1}="UL";
(3,13)*+{_2}="UR";
(0,7)*{\bu}="U";
(0,-7)*+{_{2k}}="D";
(0,-10)*{\underbrace{\ \ \ \ \ \ \ \ \ \ \ \ \ \ \ \ \ }_{2k\ \text{valent vertex}}};
    (0,0)*{\bu}="C";
  (7,0)*+{_3}="R";
   (-7,0)*+{_0}="L";
\ar @{<-} "U";"UL" <0pt>   
\ar @{<-} "U";"UR" <0pt>  
\ar @{<-} "C";"U" <0pt>
\ar @{<-} "C";"D" <0pt>
\ar @{->} "C";"L" <0pt>
\ar @{->} "C";"R" <0pt>
 \endxy}
 \Ea\hspace{-2mm}.
 $$
 Let us call such representatives of equivalence classes of generators of $\cE_0(\widehat{\Delta}\RGC_{1})$ {\it good}\, ones.
The  point is that the differential $\delta_0$ acting on a good ribbon quiver  does {\it not}\, create ribbon quivers having  edges of the above types.

\sip

Given any good generator $\Ga\in \cE_0(\widehat{\Delta}\RGC_{1})$, let $\widehat{\Ga}$ be the ribbon quiver obtained from $\Ga$ by contracting all its operadic edges. As no operadic edges in $\Ga$ are connected to the $\caD\Lie_\infty$-type vertices of $\Ga$, the ribbon quiver  $\widetilde{\Ga}$ can be identified with an element of the prop enveloping $P\caD\Lie_\infty$ of the dg properad $\caD\Lie_\infty$ to whose elements some number of sources and targets (of valencies $\geq 2$) are attached. 
As the differential acts trivially on such sources and targets, we can assume without loss of generality 
that these sources and targets (as well as the edges attached to them) are all distinguished, say all edges attached to  targets are marked by integers $\{1,\ldots, N\}$ and all edges attached to sources  are marked by
integers  $\{\bar{1},\ldots, \bar{N}\}$. Denote the complex spanned by such marked graphs $\widetilde{\Ga}$ by $A^{marked}(N,N)$. The complex $A^{marked}(N,N)$ can be identified with the tensor product of a trivial complex (corresponding to sources and targets) with the complex $P\caD\Lie_\infty(N,N)$ whose cohomology is generated \cite{L}  by 4-valent corollas  (\ref{2: fourvalent vertex}) modulo the relation (\ref{2: double Lie relations}). Using Maschke's theorem (which allows us to forget markings at the cohomology level), we conclude that there is an isomorphism of graded vector spaces,
$$
H^\bu(\widehat{\Delta}\RGC_{1}) \simeq  \cE_0 (\sRGC_1),
$$
which in turn implies that the morphism (\ref{3: map q_1}) is an isomorphism of complexes. The Proposition is proven.
\end{proof}

To finish the proof of the Main Theorem {\ref{2: Main theorem}} it is enough to observe that the epimorphism $\pi: \ORGC_{d+1} {\rar} \sRGC_{d+1}$ is equal to the composition $q\circ p$ of quasi-isomorphisms $p$ and $q$ studied in the above two Propositions.

\bip

\bip

\def\cprime{$'$}

\end{document}